\documentclass[11pt]{article}
\usepackage[a4paper,margin=25mm]{geometry}
\usepackage{xcolor}
\usepackage{amsmath,amssymb,amsthm,mathtools,booktabs,graphicx,hyperref,enumitem}
\usepackage{tikz}
\usepackage[abbrev]{amsrefs}
\hypersetup{hidelinks,pdftitle={Upwind finite element approximation of the gradient blow-up time for a viscous Hamilton--Jacobi equation},
 pdfauthor={Yohei Fujishima and Takahiro Yamada}}
\numberwithin{equation}{section}
\allowdisplaybreaks[4]
\newtheorem{theorem}{Theorem}[section]
\newtheorem{lemma}[theorem]{Lemma}
\newtheorem{proposition}[theorem]{Proposition}
\newtheorem{corollary}[theorem]{Corollary}
\theoremstyle{definition}
\newtheorem{example}[theorem]{Example}
\theoremstyle{remark}\newtheorem{remark}[theorem]{Remark}
\newcommand{\norm}[1]{\lVert #1\rVert}
\newcommand{\R}{\mathbb R}
\newcommand{\one}{\mathbf 1}
\newcommand{\dt}{\tau_{\max}}
\newcommand{\Ih}{I_h}
\newcommand{\Rh}{R_h}
\title{Upwind finite element approximation
of the gradient blow-up time
for a viscous Hamilton--Jacobi equation}
\author{Yohei Fujishima\thanks{Corresponding author.
Email: \href{mailto:fujishima@shizuoka.ac.jp}{\texttt{fujishima@shizuoka.ac.jp}}.}
\quad\mbox{and}\quad Takahiro Yamada\thanks{Email: \href{mailto:yamada.takahiro.21@shizuoka.ac.jp}{\texttt{yamada.takahiro.21@shizuoka.ac.jp}}.}\\[0.5em]
\small Department of Mathematical and Systems Engineering\\
\small Faculty of Engineering, Shizuoka University\\
\small 3-5-1 Johoku, Hamamatsu 432-8561, Japan}
\date{}

\begin{document}
\maketitle
\begin{abstract}
We analyze an upwind finite element method with piecewise linear elements for a viscous Hamilton--Jacobi equation with a superquadratic gradient source term, subject to homogeneous Dirichlet boundary conditions. On polytopal inner approximations of smooth bounded domains in dimensions~$N=1,2,3$, we prove maximum-norm convergence of the numerical solution to the exact solution before gradient blow-up and convergence of gradient-threshold times to the blow-up time under explicit mesh and regularity assumptions. We also verify the regularity assumption and finite-time gradient blow-up for an example on the disk, to which the main theorem therefore applies. Numerical experiments illustrate this gradient blow-up.
\end{abstract}

\medskip 

\noindent 
{\small
\vspace{5pt}
}

\section{Introduction}\label{Section:1}
We study the numerical detection of gradient blow-up using an upwind
finite element approximation of the
viscous Hamilton--Jacobi equation
\begin{equation}\label{eq:1.1}
 \begin{cases}
 \partial_t u=\Delta u+|\nabla u|^p
       & \mbox{in } \Omega \times (0,T^*),\\[3pt]
 u=0
       & \mbox{on } \partial\Omega \times (0,T^*),\\[3pt]
 u(\cdot,0)=u_0
       & \mbox{in }  \Omega,
 \end{cases}
\end{equation}
where $\Omega\subset\R^N$ is a bounded domain, $N\in\{1,2,3\}$ is the
spatial dimension, $p>2$, the initial value~$u_0$ vanishes on $\partial\Omega$, and $T^*$
denotes the maximal existence time of the classical solution.
This equation arises in stochastic optimal control; see
\cite{MizogS2023}.
On smooth bounded domains, the problem admits a unique local
classical solution for each $u_0\in C^1(\overline\Omega)$ satisfying the
boundary condition; see \cite{Soupl2002a}*{Section~1}.
Since the maximum principle gives~$\norm{u(\cdot,t)}_{\infty,\Omega}\le\norm{u_0}_{\infty,\Omega}$ for~$0\le t<T^*$, the solution remains uniformly bounded. The continuation criterion therefore implies
gradient blow-up whenever
$T^*<\infty$:
\begin{equation}\label{eq:1.2}
 0<T^*<\infty,\qquad
 \limsup_{t\uparrow T^*}\norm{\nabla u(\cdot,t)}_{\infty,\Omega}=\infty.
\end{equation}
Thus, $T^*$ marks the loss of classical regularity even though the
solution itself remains bounded.
The gradient estimate in~\cite{SouplZ2006}*{Theorem~3.2(i)} gives a uniform
bound on every compact subset of~$\Omega$ up to~$T^*$, so this blow-up can
occur only on the boundary. Moreover, the maximum of~$|\nabla u(\cdot,t)|$
on the boundary itself becomes unbounded as~$t\uparrow T^*$
\cite{Soupl2002a}*{Remark~3.3(a)}.
Whether such boundary blow-up occurs depends on the initial data.
In particular, sufficiently small $C^1$ data yield global classical solutions,
whereas data satisfying a suitable largeness condition lead to finite-time gradient blow-up
\cite{Soupl2002a}*{Proposition~3.1 and Theorem~2.1}.
Appendix~\ref{Section:A} gives an explicit criterion for the latter case.
Further analytical studies of \eqref{eq:1.1} address gradient blow-up
rates and profiles \cites{PorreS2017,AttouS2020,MizogS2023a}, as well as
global continuation as viscosity solutions and the loss and recovery of
boundary conditions \cites{BarleD2004,PorreS2017a,PorreS2020,MizogS2023}.

Blow-up in semilinear heat equations has been studied
extensively since Fujita's pioneering work \cite{Fujit1966}.
On the numerical side, Nakagawa proved convergence before
blow-up and convergence of numerical blow-up times using adaptive
finite differences \cite{Nakag1975}.
Abstract sufficient conditions for convergence of approximate blow-up
times are given in
\cite{Ushij2000}. Other numerical approaches include mass-lumped
finite elements for nonlinear boundary conditions
\cite{AcostDR2002}, adaptive Euler methods
\cite{Grois2006}, and space and time adaptation
based on conditional a posteriori bounds
\cite{CangiGKM2016}.
Rescaling provides another approach: for the algorithm of Berger and Kohn
\cite{BergeK1988}, convergence of the computed blow-up time
for a semilinear heat equation is established in \cite{ChoS2026}.
In computations, the first time the numerical solution reaches a large prescribed amplitude can serve as an indicator of the blow-up time~\cite{BrunnWZ2010}. Cho~\cite{Cho2013} uses an amplitude threshold that tends to infinity as the time step tends to zero.

For gradient blow-up, however, the solution values remain bounded, so their
size cannot serve as a blow-up indicator. We therefore record, on each mesh, the
first time the $L^\infty$ norm of the numerical gradient reaches or exceeds
a prescribed threshold. We show that this time converges to~$T^*$ for a
suitable choice of threshold tending to infinity as the mesh size tends to zero.

\subsection{The numerical scheme}
We write the nonlinear source as a drift term~$|\nabla u|^p=b(\nabla u)\cdot\nabla u$, where~$b(q)=|q|^{p-2}q$. To discretize this term while preserving the maximum principle, we combine Tabata's upwind finite element construction~\cite{Tabat1977}*{Section~5} with mass lumping. For related conservative upwind methods, see Baba and Tabata~\cite{BabaT1981} and~\cite{Knobl2006}*{Sections~3--4}.
Let $\mathcal T_h$ be a simplicial mesh of a polytopal
inner approximation $\Omega_h$ of the smooth domain $\Omega$, with
$h:=\max_{K\in\mathcal T_h}\operatorname{diam}K$, and let $U_h^n$ be the
continuous piecewise affine numerical solution at time $t_n$.
For a possibly nonuniform time grid $0=t_0<t_1<\cdots$, we write
\begin{equation}\label{eq:1.3}
 \tau_n:=t_{n+1}-t_n>0,\qquad
 \dt:=\sup_{n\ge0}\tau_n.
\end{equation}
At each interior vertex, the averaged gradient of~$U_h^n$ determines the drift and hence the upwind element. Using the drift and upwind elements determined at time~$t_n$, we formulate a linear system for~$U_h^{n+1}$. Under a suitable mesh condition, this system satisfies a discrete maximum principle (DMP).

To track the growth of the numerical gradient, we use its maximum over the mesh elements,
\begin{equation*}
 G_h^n:=\norm{\nabla U_h^n}_{\infty,\Omega_h}
       =\max_{K\in\mathcal T_h}\bigl|(\nabla U_h^n)_K\bigr|.
\end{equation*}
On the shape-regular, quasiuniform meshes considered here, the DMP and an inverse estimate bound~$hG_h^n$ uniformly in~$h$ and~$n$. Thus~$G_h^n$ remains bounded on each fixed mesh, although it may grow as the mesh is refined. A fixed threshold may be crossed well before~$T^*$, so we let the threshold grow with refinement. For a positive scale~$L_h\to\infty$ as~$h\downarrow0$ and a multiplier~$\kappa>0$, we define the first crossing time
\begin{equation*}
 T_h(\kappa):=\inf\{t_n\colon G_h^n\ge\kappa L_h\}.
\end{equation*}
Our aim is to find conditions relating~$L_h$, $h$, and~$\dt$ under which the gradient-threshold times~$T_h(\kappa)$ converge to the blow-up time~$T^*$ uniformly for~$\kappa$ in any fixed compact subinterval of~$(0,\infty)$.

\subsection{Outline of the main results}
To determine when the numerical gradient first reaches the threshold,
we begin with estimates on each fixed interval~$[0,T]$ with~$T<T^*$. 
Under the hypotheses of Theorem~\ref{Theorem:3.5}, the condition~$\dt h^{-1}\to0$ as~$h\downarrow0$ yields
\begin{equation*}
 \max_{t_n\le T}\norm{U_h^n-u(t_n)}_{\infty,\Omega_h}\le C_T(h+\dt)\to0,
 \qquad \max_{t_n\le T}G_h^n\le C_T,
\end{equation*}
where~$C_T$ is independent of~$h$.
By the gradient bound, any diverging threshold stays above~$G_h^n$ throughout~$[0,T]$ for all sufficiently small~$h$. On the other hand, if gradient blow-up occurs as in~\eqref{eq:1.2}, uniform convergence of the numerical solution implies that~$G_h^n$ exceeds any fixed level before~$T^*$ on all sufficiently fine meshes. In Theorem~\ref{Theorem:4.2}, we combine these facts to show that a threshold scale~$L_h\to\infty$ can be chosen so that
\begin{equation}\label{eq:1.4}
 \sup_{\kappa\in[\kappa_-,\kappa_+]}|T_h(\kappa)-T^*|\to0
 \quad\text{for every }0<\kappa_-\le\kappa_+<\infty.
\end{equation}
Theorem~\ref{Theorem:4.2} takes the uniform convergence of the solution values and the numerical gradient bound on each fixed interval as hypotheses. Beyond these, it applies on bounded Lipschitz domains and requires only that~$u$ and~$\nabla u$ be continuous on~$\overline\Omega\times[0,T]$ for each~$T<T^*$. However, the construction gives no explicit rate at which~$L_h$ diverges as~$h\downarrow0$.

In Theorem~\ref{Theorem:4.6}, assuming that the exact solution satisfies the regularity condition~\hyperlink{regularity-R}{(R)} of Section~\ref{Subsection:3.2}, we show that the uniform convergence in~\eqref{eq:1.4} holds whenever the threshold scale~$L_h$ satisfies
\begin{equation*}
 L_h\to\infty,\qquad
 hL_h^{3(p-1)}\to0,\qquad
 \dt h^{-1}L_h^{p-1}\to0 \quad\text{as }h\downarrow0.
\end{equation*}
For example, when~$\dt=h^2$, these conditions allow~$L_h=h^{-\sigma}$ for any~$0<\sigma<1/(3(p-1))$.
Thus, in contrast to Theorem~\ref{Theorem:4.2}, Theorem~\ref{Theorem:4.6} gives an explicit choice of~$L_h$ for which the uniform convergence in~\eqref{eq:1.4} holds.
Moreover, Proposition~\ref{Proposition:3.2} shows that, on a smooth bounded domain, suitable conditions on the initial data ensure~\hyperlink{regularity-R}{(R)}.

\subsection{Strategy}
Both convergence proofs split into a lower bound and an upper bound for~$T_h(\kappa)$, and the lower bound is the same in both. 
For each fixed~$T<T^*$, the numerical gradient is uniformly bounded on~$[0,T]$, so~$T_h(\kappa)>T$ for all sufficiently small~$h$ whenever~$L_h\to\infty$ and~$\kappa>0$ is fixed. For the upper bound in Theorem~\ref{Theorem:4.2}, we compare the spatial difference quotients of the exact and numerical solutions near the blow-up time~$T^*$.
Where the exact gradient is large, we can find a short segment along which the exact solution has a large difference quotient. With the endpoints fixed, the numerical values at the nearest preceding time level converge to the exact values at both endpoints. The numerical difference quotient is therefore large as well, so the numerical gradient must be large on at least one mesh element. Applying this comparison at a sequence of times tending to~$T^*$ gives, by a diagonal argument, a scale~$L_h\to\infty$ such that, for each fixed~$\kappa>0$, the threshold~$\kappa L_h$ is crossed before~$T^*$ on all sufficiently fine meshes.

When~$L_h$ is prescribed, as in Theorems~\ref{Theorem:4.5} and~\ref{Theorem:4.6}, we must show that the numerical gradient reaches the specified level~$\kappa L_h$ before~$T^*$. We use the same comparison on a segment chosen short enough that the exact directional derivative stays large along it. As the segment shrinks with~$h$, the endpoint errors must be small relative to the change in the exact solution along the segment (see~\eqref{eq:4.12}). To control these errors, we use a maximum-norm estimate on~$\Omega_h$ that accounts for both the discretization error and the discrepancy at the computational boundary. Following the projection approach of Thom\'ee~\cite{Thome2006}*{Chapter~1}, we extend the exact solution to a fixed outer domain, take its Ritz projection there, and compare the projected solution with~$U_h^n$ on~$\Omega_h$ at time~$t_n$. Under the conditions on~$h$, $\dt$, and~$L_h$, we obtain a quantitative maximum-norm bound on this error up to the first threshold crossing (see~\eqref{eq:3.26}), and we use this bound to show that the endpoint errors are negligible relative to the change along the segment. The numerical difference quotient then forces a crossing before~$T^*$. Together with the lower bound, this gives~$T_h(\kappa)\to T^*$.

\paragraph{Organization of the paper.}
Section~\ref{Section:2} presents the scheme on $\Omega_h$, the geometric
and mesh conditions, nonlinear error stability, Poincar\'e inequalities,
and approximation estimates. Section~\ref{Section:3} then establishes a reference-function
error bound and value convergence on inner approximations.
Building on these estimates, Section~\ref{Section:4} proves that there exists a threshold scale for which the first crossing times converge to~$T^*$ and gives explicit sufficient conditions for this convergence when the scale is prescribed.
Section~\ref{Section:5} describes the disk meshes and presents the numerical experiments.
Appendix~\ref{Section:A} establishes regularity estimates and a large-data blow-up criterion.

\paragraph{Notation.}
We write $B(x,r)$ for the open Euclidean ball centered at $x$ with radius $r$.
For $1\le q\le\infty$ and a measurable set $E$, write $\norm{f}_{q,E}:=\norm{f}_{L^q(E)}$. For a mesh element $K$, we write $(f)_K:=f|_K$ for the restriction of $f$ to $K$. When this restriction is constant, we identify it with its constant value; in particular, $(\nabla v_h)_K$ is a constant vector when $v_h$ is piecewise affine.
The symbol $|\cdot|$ denotes the Euclidean norm for vectors and the Frobenius norm for matrices and higher-order tensors. For a finite vector $v=(v_r)$, also write $|v|_1=\sum_r|v_r|$. For a matrix~$X$, $\norm{X}_{r\to s}$ denotes its operator norm from~$\ell^r$ to~$\ell^s$.
Matrix and vector inequalities are entrywise. Constants $c,C>0$ may change from line to line and are independent of $h$ and $T<T^*$ unless otherwise indicated. For nonnegative quantities $a,b$, we write $a\lesssim b$ (or $b\gtrsim a$) if $a\le Cb$ with such a constant $C$, and $a\simeq b$ if both $a\lesssim b$ and $b\lesssim a$ hold.
For fixed $T<T^*$, we write $O_T(h)$ for a quantity bounded in absolute value by $C_Th$, with $C_T$ independent of $h$.

\section{The scheme and preliminary estimates}\label{Section:2}
We construct the upwind finite element scheme on the computational domain
$\Omega_h$ and prove its discrete maximum principle. To control errors
between two solutions, we then relate changes in the nonlinear drift to
the diffusion coefficients. This gives a condition under which each time
step is contractive in the maximum norm. The final subsection collects
the approximation estimates needed to apply this stability result to the
continuous solution.

\subsection{The upwind scheme and discrete maximum principle}
Let $\mathcal T_h$ be a conforming simplicial mesh of a polytopal
computational domain $\Omega_h\subset\Omega$. In the main setting,
these domains approximate the fixed smooth domain $\Omega$ from inside
and satisfy the following geometric condition:
\begin{equation}\label{eq:2.1}
 \Omega_h\subset\Omega,\qquad
 \sup_{x\in\Omega\setminus\Omega_h}
       \operatorname{dist}(x,\partial\Omega)\lesssim h.
\end{equation}
For classical domain approximations, see
\cite{Tabat1978}*{Section~1 and Lemma~3.1} and
\cite{Thome2006}*{Chapter~1}.
On this mesh, let $\mathcal J_h=\mathcal I_h\cup\mathcal B_h$ index the interior and
boundary vertices of $\Omega_h$, denoted by $x_i$.
With $P_1(K)$ denoting the affine functions on $K$, define
\begin{equation*}
 V_h=\{v_h\in C(\overline\Omega_h)\colon (v_h)_K\in P_1(K)
                  \text{ for all }K\in\mathcal T_h\},\qquad
 V_h^0=V_h\cap H_0^1(\Omega_h).
\end{equation*}
The corresponding nodal basis $\{\varphi_i\colon i\in\mathcal J_h\}$ satisfies
$\varphi_i(x_j)=\delta_{ij}$, $\varphi_i\ge0$, and $\sum_i\varphi_i=1$.
Thus, a nodal vector $W=(W_i)_{i\in\mathcal J_h}$ represents
\begin{equation}\label{eq:2.2}
 W_h=\sum_{i\in\mathcal J_h}W_i\varphi_i\in V_h,\qquad
 \norm{W_h}_{\infty,\Omega_h}=\norm{W}_{\ell^\infty}
                  =\max_{i\in\mathcal J_h}|W_i|.
\end{equation}
We identify finite element functions with their nodal vectors when
applying the discrete operators defined below.
Using this basis, we also define the \emph{nodal interpolation} of any
continuous function $f$ on $\overline\Omega_h$ by
\begin{equation}\label{eq:2.3}
 \Ih:C(\overline\Omega_h)\to V_h,\qquad
 \Ih f=\sum_{i\in\mathcal J_h}f(x_i)\varphi_i.
\end{equation}
If $f$ vanishes on $\partial\Omega_h$, then $\Ih f\in V_h^0$.
When interpolating $H^2$ functions, we use the continuous
representative, which exists for $N=1,2,3$.

To define the discrete diffusion, we first introduce the lumped masses and stiffness matrix entries
\begin{equation}\label{eq:2.4}
 m_i=\int_{\Omega_h}\varphi_i\, \mathrm{d}x>0,\qquad
 A_{ij}=\int_{\Omega_h}\nabla\varphi_i\cdot\nabla\varphi_j\, \mathrm{d}x,
 \qquad \sum_{j\in\mathcal J_h}A_{ij}=0.
\end{equation}
For an interior vertex $i$, let $\omega_i=\bigcup_{K\ni x_i}K$ be its vertex patch and
let $\mathcal N_i$ index its edge neighbors.
In terms of this patch, the lumped mass is given by
\begin{equation*}
 m_i=\frac{1}{N+1}\sum_{K\ni x_i}|K|=\frac{|\omega_i|}{N+1}.
\end{equation*}
Here $|K|$ and $|\omega_i|$ denote the $N$-dimensional volumes of
$K$ and $\omega_i$, respectively.
Set $a_{ij}=-A_{ij}/m_i$ for $i\in\mathcal I_h$ and $j\ne i$,
and $a_{ij}=0$ for $i\in\mathcal B_h$ and $j\ne i$. Following~\cite{Thome2006}*{Chapter~15}, we define the corresponding lumped discrete Laplacian for an arbitrary nodal vector $V$ by
\begin{equation}\label{eq:2.5}
 (\Delta_h V)_i=\begin{cases}
 -\displaystyle \frac{1}{m_i}\sum\limits_{j\in\mathcal J_h}A_{ij}V_j
       =\sum\limits_{j\in\mathcal N_i}a_{ij}(V_j-V_i),&i\in\mathcal I_h,\\
 0,&i\in\mathcal B_h.
 \end{cases}
\end{equation}
For interior nodes, the second equality in~\eqref{eq:2.5} follows from the definition of $a_{ij}$, the zero row sum in~\eqref{eq:2.4}, and the fact that $A_{ij}=0$ whenever $j\notin\mathcal N_i\cup\{i\}$.
We identify $\Delta_h V$ with its finite element
function, so that $\Delta_h:V_h\to V_h^0$.

To discretize the drift term, recall that $b(q)=|q|^{p-2}q$ for $q\in\R^N$,
so that $b(q)\cdot q=|q|^p$. We then define the recovered gradient
$g_i(W)$ and the associated drift $b_i(W)$ at interior nodes by
\begin{equation}\label{eq:2.6}
 g_i(W)=\frac{1}{|\omega_i|}\sum_{K\ni x_i}|K|(\nabla W_h)_K,
 \qquad b_i(W)=b(g_i(W)).
\end{equation}
The weights $|K|/|\omega_i|$ are nonnegative and sum
to one. Consequently, for every nodal vector $W$ and
$v\in C^1(\overline\Omega_h)$,
\begin{equation}\label{eq:2.7}
 |g_i(W)-\nabla v(x_i)|
 \le\max_{K\ni x_i}|(\nabla W_h)_K-\nabla v(x_i)|
 \le\norm{\nabla(W_h-v)}_{\infty,\Omega_h}.
\end{equation}
For each $i\in\mathcal I_h$, choose a simplex $K_i(b_i(W))$ in the patch containing
$x_i+s b_i(W)$ for small $s>0$.
If the ray lies on a shared face, continuity of finite element functions
ensures that the directional derivative along the ray is the same
in each simplex sharing that face, so any of them may be chosen.
If $b_i(W)=\mathbf0$, any simplex in the patch may be chosen.
For fixed $W$, we now express the drift in matrix form following
\cite{Tabat1978}.
Let $x_i,x_{j_1},\ldots,x_{j_N}$
be the vertices of the selected simplex $K=K_i(b_i(W))$ and define
$E_{i,K}$ by taking its edge vectors as columns. With this notation,
for $v_h\in V_h$ with nodal vector $V=(V_j)_{j\in\mathcal J_h}$, where $V_j=v_h(x_j)$, we have
\begin{equation}\label{eq:2.8}
 E_{i,K}=\bigl(x_{j_1}-x_i,\ldots,x_{j_N}-x_i\bigr),\qquad
 E_{i,K}^{\mathsf{T}}(\nabla v_h)_K=(V_{j_1}-V_i,\ldots,V_{j_N}-V_i)^{\mathsf{T}}.
\end{equation}
By the choice of $K$,
\begin{equation}\label{eq:2.9}
 \gamma:=E_{i,K}^{-1}b_i(W)\ge\mathbf0.
\end{equation}
Set $\gamma_{i,j_r}(W):=\gamma_r$ for $1\le r\le N$ and the other entries to zero.
Then, by \eqref{eq:2.8} and \eqref{eq:2.9},
\begin{equation}
  \label{eq:2.10}
 \begin{aligned}
 b_i(W)&=\sum_{j\in\mathcal N_i}\gamma_{ij}(W)(x_j-x_i),\\
 (B_h(W)V)_i&:=b_i(W)\cdot(\nabla v_h)_{K_i(b_i(W))}
       =\sum_{j\in\mathcal N_i}\gamma_{ij}(W)(V_j-V_i).
 \end{aligned}
\end{equation}
Finally, for boundary nodes $i\in\mathcal B_h$, set $(B_h(W)V)_i=0$
and $\gamma_{ij}(W)=0$.

Combining the diffusion and drift operators, we define the matrix on all nodal vectors by
\begin{equation}\label{eq:2.11}
 Q(W):=\Delta_h+B_h(W).
\end{equation}
Writing $Q(W)=(q_{ij}(W))$ and using~\eqref{eq:2.5},~\eqref{eq:2.10}, and~\eqref{eq:2.11}, we obtain its off-diagonal and diagonal entries as follows:
\begin{equation}\label{eq:2.12}
 q_{ij}(W)=a_{ij}+\gamma_{ij}(W),\qquad
 q_{ii}(W)=-\sum_{j\ne i}q_{ij}(W).
\end{equation}
Here the first formula holds for~$i\ne j$. In particular, $Q(W)\one=\mathbf0$, where~$\one$ is the nodal vector with all entries equal to one.

We next introduce the time discretization. For each $h$, assume that the time grid in \eqref{eq:1.3}
has only finitely many nodes in every bounded time interval.
On this grid, let $U_h^n\in V_h^0$ denote the numerical solution at time $t_n$,
with nodal vector $U^n=(U_i^n)_{i\in\mathcal J_h}$, where
$U_i^n=U_h^n(x_i)$. With this notation, the scheme is
\begin{equation}\label{eq:2.13}
 m_i\frac{U_i^{n+1}-U_i^n}{\tau_n}+\sum_{j\in\mathcal J_h}A_{ij}U_j^{n+1}
   -m_i b_i(U^n)\cdot(\nabla U_h^{n+1})_{K_i(b_i(U^n))}=0
\end{equation}
for $i\in\mathcal I_h$, with $U_i^n=U_i^{n+1}=0$ for $i\in\mathcal B_h$.
Since the old vector $U^n$ fixes the drift and upwind elements, each step is linear in $U^{n+1}$.
Equivalently, dividing~\eqref{eq:2.13} by~$m_i$ and using~\eqref{eq:2.5},~\eqref{eq:2.10}, and~\eqref{eq:2.11} gives the matrix form
\begin{equation}\label{eq:2.14}
 (I-\tau_nQ(U^n))U^{n+1}=U^n.
\end{equation}
To initialize the scheme, set $U_i^0=u_0(x_i)$ for $i\in\mathcal I_h$ and $U_i^0=0$
for $i\in\mathcal B_h$.
For an inner approximation, the discrepancy between $u_0$ and the zero
computational boundary values is included in the error estimates of
Section~\ref{Subsection:3.2}.

The following theorem establishes a discrete maximum principle for the scheme~\eqref{eq:2.13}. To this end, we assume the mesh sign condition
\begin{equation}\label{eq:2.15}
 a_{ij}\ge0\quad\text{for }i\in\mathcal I_h,\quad j\in\mathcal J_h\setminus\{i\}.
\end{equation}
Equivalently, $A_{ij}\le0$ for the same indices.
The proof uses matrix positivity and maximum-norm contractivity; see
\cite{Thome2006}*{Theorems~15.5--15.6} for the lumped-mass heat equation
and \cite{ThomeW2008}*{Theorems~3.1--3.3} for related matrix criteria.

\begin{theorem}[Discrete maximum principle]\label{Theorem:2.1}
Assume the sign condition~\eqref{eq:2.15}.
Then, for every $U_h^n\in V_h^0$ and $\tau_n>0$, the scheme
\eqref{eq:2.13} has a unique solution $U_h^{n+1}\in V_h^0$, which satisfies
\begin{equation*}
 \min\{0,\min_jU_j^n\}\le U_h^{n+1}(x)\le\max\{0,\max_jU_j^n\}
\end{equation*}
for $x\in\overline\Omega_h$.
\end{theorem}
\begin{proof}
Fix $U^n$ and $\tau=\tau_n>0$.
By \eqref{eq:2.9}, $\gamma_{ij}(U^n)\ge0$ for $i\ne j$.
Together with~\eqref{eq:2.15}, \eqref{eq:2.12} therefore implies that the
off-diagonal entries $q_{ij}(U^n)$ are nonnegative.
Using \eqref{eq:2.12}, we obtain, for every nodal vector $V$ and every $i\in\mathcal J_h$,
\begin{equation}\label{eq:2.16}
 \bigl((I-\tau Q(U^n))V\bigr)_i=V_i+\tau\sum_{j\ne i}q_{ij}(U^n)(V_i-V_j).
\end{equation}
Then, choosing $i$ with $|V_i|=\norm{V}_{\ell^\infty}$ and using $q_{ij}(U^n)\ge0$ for $j\ne i$, we obtain
\[
 \bigl|\bigl((I-\tau Q(U^n))V\bigr)_i\bigr|
 \ge |V_i|+\tau\sum_{j\ne i}q_{ij}(U^n)(|V_i|-|V_j|)
 \ge\norm{V}_{\ell^\infty}.
\]
Thus $I-\tau Q(U^n)$ is invertible and
\begin{equation}\label{eq:2.17}
 \norm{(I-\tau Q(U^n))^{-1}V}_{\ell^\infty}\le\norm{V}_{\ell^\infty}
 \quad\text{for every nodal vector }V.
\end{equation}
For $W\ge\mathbf0$, set $V=(I-\tau Q(U^n))^{-1}W$.
A negative minimum $V_i=\min_jV_j<0$ would give $W_i\le V_i<0$ by~\eqref{eq:2.16}, a contradiction.
Hence $(I-\tau Q(U^n))^{-1}\ge O$.
Since $Q(U^n)\one=\mathbf0$, we also have $(I-\tau Q(U^n))^{-1}\one=\one$.
Each component of the unique solution $U^{n+1}=(I-\tau Q(U^n))^{-1}U^n$ is therefore a convex combination of the components of $U^n$.
Hence
\[
 \min_jU_j^n\le U_i^{n+1}\le\max_jU_j^n,
 \qquad i\in\mathcal J_h.
\]
The zero boundary rows of $Q(U^n)$ in \eqref{eq:2.11} give $U_i^{n+1}=U_i^n=0$ for $i\in\mathcal B_h$.
Finally, $\varphi_i\ge0$ and $\sum_i\varphi_i=1$ extend the nodal bounds to $U_h^{n+1}$ throughout $\overline\Omega_h$.
\end{proof}

\subsection{Mesh conditions and nonlinear error stability}\label{Subsection:2.2}
To estimate the error, we must control how the drift varies with the
solution. We express this nonlinear difference as a matrix acting
on the nodal error and give mesh and time-step conditions that make the
resulting error propagation contractive in the maximum norm.

Put $h_K=\operatorname{diam}K$, $h=\max_{K\in\mathcal T_h}h_K$, and let $r_K$ be the inradius. In this paper, we assume shape regularity and quasiuniformity in the form
\begin{equation}\label{eq:2.18}
 \frac{h_K}{r_K}\le C,\qquad ch\le h_K\le h, 
\end{equation}
for every element on every mesh.
The mesh condition~\eqref{eq:2.18} gives, for~$E_{i,K}$ from~\eqref{eq:2.8}, the following bounds (see, for example,~\cite{BrennS2008}*{Section~4.4}):
\begin{equation}\label{eq:2.19}
 c h_K^N\le |K|\le C h_K^N,\qquad
 \norm{E_{i,K}}_{2\to2}\le C h_K,\qquad \norm{E_{i,K}^{-1}}_{2\to2}\le C h_K^{-1}\le C h^{-1}.
\end{equation}
In particular, the third bound gives $|(\nabla\varphi_j)_K|\lesssim h^{-1}$ for every vertex $x_j$ of $K$.
Patch simplices have disjoint interiors, lie in $B(x_i,h)$,
and each have volume at least $ch^N$, so their number is uniformly bounded.
Consequently, with constants also depending on~$N$,
\begin{equation}\label{eq:2.20}
 m_i\simeq h^N.
\end{equation}
To express the recovered gradients $g_i(W)$ in terms of nodal values, we insert \eqref{eq:2.2} into
\eqref{eq:2.6} and interchange the two finite sums. This gives
\begin{equation*}
 g_i(W)=\sum_{j\in\mathcal N_i\cup\{i\}}\ell_{ij}W_j,
 \qquad \ell_{ij}=\frac{1}{|\omega_i|}\sum_{K\ni x_i,x_j}|K|(\nabla\varphi_j)_K\in\R^N.
\end{equation*}
Set $\ell_{ij}=\mathbf0$ for other $j$ and for boundary rows.
The second and third identities below show, respectively, that gradient recovery gives zero for constant nodal data and that it reproduces the gradient of every affine function. We use the vanishing diagonal later in the formula for the nonlinear difference.
\begin{lemma}\label{Lemma:2.2}
For every interior patch,
\begin{equation}\label{eq:2.21}
 \ell_{ii}=\mathbf0,\qquad \sum_j\ell_{ij}=\mathbf0,\qquad
 \sum_j\ell_{ij}(x_j-x_i)^{\mathsf{T}}=I_N.
\end{equation}
Under \eqref{eq:2.18}, $\sum_j|\ell_{ij}|\lesssim h^{-1}$.
\end{lemma}
\begin{proof}
Since $\varphi_i$ vanishes on $\partial\omega_i$,
the divergence theorem gives
\[
 |\omega_i|\ell_{ii}=\int_{\omega_i}\nabla\varphi_i\, \mathrm{d}x=\mathbf0.
\]
This proves the first identity in \eqref{eq:2.21}.
Taking $W=\one$ and then $W_j=\xi\cdot(x_j-x_i)$ in \eqref{eq:2.6},
with arbitrary $\xi\in\R^N$, gives $g_i(W)=\mathbf0$ and $g_i(W)=\xi$, respectively.
The representation of $g_i(W)$ above then yields the remaining identities in
\eqref{eq:2.21}. Finally, \eqref{eq:2.19} yields
\[
 \sum_j|\ell_{ij}|\le \frac{1}{|\omega_i|}
 \sum_{K\ni x_i}|K|\sum_{x_j\in K}\bigl|(\nabla\varphi_j)_K\bigr|
 \lesssim \frac{h^{-1}}{|\omega_i|}\sum_{K\ni x_i}|K|=h^{-1}.
 \qedhere
\]
\end{proof}

The next bounds control the row sums of the diffusion and upwind
coefficients in terms of the mesh size and the norm of the drift.
\begin{proposition}\label{Proposition:2.3}
Assume the sign condition~\eqref{eq:2.15} and the mesh condition~\eqref{eq:2.18}. Then the coefficients satisfy
\begin{equation*}
 \sum_{j\ne i}a_{ij}\lesssim h^{-2}\quad\mbox{for}\ i\in\mathcal J_h,\qquad
 \sum_j\gamma_{ij}(W)\lesssim h^{-1}|b_i(W)|\quad\mbox{for}\ i\in\mathcal I_h.
\end{equation*}
\end{proposition}
\begin{proof}
For $i\in\mathcal I_h$, \eqref{eq:2.4} and \eqref{eq:2.19} give the first bound.
For $K=K_i(b_i(W))$, let $\gamma$ be defined by \eqref{eq:2.9}.
Its nonnegativity, the definition of $\gamma_{ij}(W)$, and the bound~\eqref{eq:2.19} yield
\[
 \sum_j\gamma_{ij}(W)=|\gamma|_1
 \le\sqrt N\,\norm{E_{i,K}^{-1}}_{2\to2}|b_i(W)|
 \lesssim h^{-1}|b_i(W)|.
\]
For $i\in\mathcal B_h$, the first bound follows from $a_{ij}=0$ for $j\ne i$.
\end{proof}

To relate gradient recovery to diffusion, we impose, in addition to the
sign condition~\eqref{eq:2.15}, the following condition,
with $c_D$ independent of $h$:
\begin{equation}\label{eq:2.22}
 |\ell_{ij}|\le c_D h\,a_{ij}
\end{equation}
for $i\in\mathcal I_h$ and $j\ne i$.
A sufficient condition for \eqref{eq:2.22} is
\begin{equation}\label{eq:2.23}
 -A_{ij}\ge a_*h^{N-2}>0
\end{equation}
for $i\in\mathcal I_h$ and $j\in\mathcal N_i$, with $a_*$ fixed.
Indeed, the definition $a_{ij}=-A_{ij}/m_i$, together with~\eqref{eq:2.20} and~\eqref{eq:2.23}, gives $a_{ij}\gtrsim h^{-2}$, whereas
Lemma~\ref{Lemma:2.2} gives $|\ell_{ij}|\lesssim h^{-1}$.
The following examples describe mesh families satisfying these conditions.

\begin{example}[Interval meshes]
  \label{Example:2.4}
For $N=1$, every quasiuniform interval mesh, as defined by the second condition in~\eqref{eq:2.18}, is automatically shape regular. That is, the first condition in~\eqref{eq:2.18} is satisfied with $h_K/r_K=2$. It also satisfies the sign condition~\eqref{eq:2.15} and condition~\eqref{eq:2.22}.
Indeed, order the nodes increasingly and
write $\lambda_-=x_i-x_{i-1}$, $\lambda_+=x_{i+1}-x_i$ at an interior node. Then
\begin{equation*}
 \begin{gathered}
 m_i=\frac{\lambda_-+\lambda_+}{2},\quad
 a_{i,i\pm1}=\frac{2}{\lambda_\pm(\lambda_-+\lambda_+)},\quad
 \ell_{i,i-1}=-\frac{1}{\lambda_-+\lambda_+},\quad
 \ell_{i,i+1}=\frac{1}{\lambda_-+\lambda_+},\quad \ell_{ii}=\mathbf0.
 \end{gathered}
\end{equation*}
These formulas give
$|\ell_{i,i\pm1}|/a_{i,i\pm1}=\lambda_\pm/2\le h/2$, so
\eqref{eq:2.22} holds with $c_D=1/2$.
\end{example}

\begin{example}[Uniformly acute meshes]
  \label{Example:2.5}
For $N=2,3$, every family of uniformly acute meshes satisfying the shape-regularity and quasiuniformity conditions~\eqref{eq:2.18} satisfies~\eqref{eq:2.23}. Indeed, for the distances $H_i,H_j$ from the vertices $x_i,x_j$ to the hyperplanes containing their opposite facets and the interior angle $\theta_{ij}$ between those facets,
\begin{equation*}
 \int_K\nabla\varphi_i\cdot\nabla\varphi_j\,\mathrm{d}x
       =-\frac{|K|\cos\theta_{ij}}{H_iH_j}.
\end{equation*}
This identity, the mesh condition~\eqref{eq:2.18}, the bounds~\eqref{eq:2.19}, and $\theta_{ij}\le\pi/2-\varepsilon$ give the required negative contribution of order~$h^{N-2}$.
Such families also exist in three dimensions: the acute cube tetrahedralization of~\cite{VandeHZG2010}*{Sections~2 and~5}, scaled into $J^3$ subcubes and reflected across shared faces to match their triangulations, yields a uniformly acute, shape-regular, quasiuniform family on the cube with $h\simeq J^{-1}$.
\end{example}

Changing the drift can change the selected upwind simplex. To account
for this change, let $Z$ be a nodal vector and $i$ an interior node, and define
$D_i^Z(b)=b\cdot(\nabla Z_h)_{K_i(b)}$ for $b\ne\mathbf0$, and $D_i^Z(\mathbf0)=0$.
This function is linear on each cone generated by a simplex at the vertex
and continuous across its interfaces, since the restrictions of $Z_h$ from
adjacent elements agree on their common face.
Integration along the segment joining $b$ and $c$, split at the cone
interfaces, therefore gives
\begin{equation}\label{eq:2.24}
 |D_i^Z(b)-D_i^Z(c)|\le\norm{\nabla Z_h}_{\infty,\Omega_h}|b-c|.
\end{equation}
This Lipschitz bound leads to the following exact representation of
the nonlinear difference.
\begin{lemma}\label{Lemma:2.6}
Let $W,Y,Z$ be nodal vectors with gradient bounds $M_1,M_2\ge0$ such that
\begin{equation}\label{eq:2.25}
 \max\{\norm{\nabla W_h}_{\infty,\Omega_h},\norm{\nabla Y_h}_{\infty,\Omega_h}\}\le M_1,
 \qquad \norm{\nabla Z_h}_{\infty,\Omega_h}\le M_2.
\end{equation}
There exists a matrix $\mathcal C=\mathcal C(W,Y;Z)$, indexed by all nodes, such that
\begin{equation*}
 \bigl(Q(W)-Q(Y)\bigr)Z=\mathcal C(W-Y),\quad \mathcal C_{ii}=0,\quad \mathcal C\one=\mathbf0,
 \quad |\mathcal C_{ij}|\le(p-1)M_1^{p-2}M_2|\ell_{ij}|.
\end{equation*}
Its boundary rows vanish, and $\mathcal C_{ij}=0$ whenever $i\in\mathcal I_h$ and $j\notin\mathcal N_i$.
\end{lemma}
\begin{proof}
By \eqref{eq:2.11}, the diffusion terms cancel in the difference.
Thus, for every $i\in\mathcal I_h$,
\eqref{eq:2.10} gives
\begin{equation*}
 \bigl((Q(W)-Q(Y))Z\bigr)_i
 =D_i^Z(b(g_i(W)))-D_i^Z(b(g_i(Y))).
\end{equation*}
Put $W(s)=Y+s(W-Y)$ and $g_i(s)=g_i(W(s))$ for $0\le s\le1$.
Then $|g_i(s)|\le M_1$ and
$g_i'(s)=\sum_j\ell_{ij}(W_j-Y_j)$.
To apply the chain rule, note that $b$ is continuously differentiable for~$p>2$, with~$Db(\mathbf0)=O$ and
\begin{equation*}
 Db(q)=|q|^{p-2}I_N+(p-2)|q|^{p-4}qq^{\mathsf{T}}
\end{equation*}
for $q\in\R^N\setminus\{\mathbf0\}$. Hence
\begin{equation}\label{eq:2.26}
 \norm{Db(q)}_{2\to2}\le(p-1)|q|^{p-2}
\end{equation}
for every $q\in\R^N$.
The affine dependence of~$g_i(s)$ on~$s$, together with~\eqref{eq:2.24} and~\eqref{eq:2.26}, shows that the function $s\mapsto D_i^Z(b(g_i(s)))$ is Lipschitz continuous, hence absolutely continuous.
We divide $[0,1]$ into subintervals on which the upwind simplex can
be kept fixed, so that the corresponding element gradient is
a constant vector on each such subinterval.
Indeed, since $g_i(s)$ is affine and each cone is convex, the set of $s$
for which $g_i(s)$ lies in a given cone is an interval; there are only
finitely many cones, so a finite subdivision suffices.
Moreover, $b(g)$ is a nonnegative multiple of $g$, so $b(g_i(s))$ lies
in the same cone as $g_i(s)$.
Thus, on each subinterval we may choose $K_i(b(g_i(s)))$ independently
of $s$, and the chain rule gives
\begin{equation*}
 \frac{d}{ds}D_i^Z(b(g_i(s)))
 =(\nabla Z_h)_{K_i(b(g_i(s)))}\cdot Db(g_i(s))\sum_j\ell_{ij}(W_j-Y_j).
\end{equation*}
If $b(g_i(s))$ lies in a common face of the cones throughout a subinterval, 
the right-hand side is independent of
which element is chosen, since the tangential
derivatives of $Z_h$ agree across common faces.
The identity also holds at $g_i(s)=\mathbf0$ since $Db(\mathbf0)=O$.
We may therefore integrate over $0\le s\le1$ to obtain the claimed representation by setting
\begin{equation*}
 \mathcal C_{ij}=\int_0^1(\nabla Z_h)_{K_i(b(g_i(s)))}\cdot Db(g_i(s))\ell_{ij}\,\mathrm{d}s
\end{equation*}
for~$i\in\mathcal I_h$ and~$j\in\mathcal J_h$, and taking the boundary rows to be zero. Finally, the coefficient bound follows from \eqref{eq:2.25}--\eqref{eq:2.26}
and this formula, while \eqref{eq:2.21} gives $\mathcal C_{ii}=0$ and $\mathcal C\one=\mathbf0$.
\end{proof}

The next lemma gives a sufficient condition for the propagation of nodal errors to be contractive in the maximum norm.
\begin{lemma}\label{Lemma:2.7}
Let $Q=(q_{ij})$ satisfy $Q\one=\mathbf0$ and $q_{ij}\ge0$ for $i\ne j$. Suppose $\mathcal C\one=\mathbf0$, $\mathcal C_{ii}=0$, and $\mathcal C_{ij}=0$ whenever $i\ne j$ and $q_{ij}=0$. For a step $\tau>0$ define
\begin{equation}\label{eq:2.27}
 \delta_Q=\tau\max_i\sum_{j\ne i}q_{ij},\qquad
 \zeta=\max_{i\ne j,\,q_{ij}>0}\frac{(-\mathcal C_{ij})_+}{q_{ij}}.
\end{equation}
Here we use the convention $\max\varnothing=0$ in the definition of $\zeta$.
If $\zeta(1+\delta_Q)\le1$, then
\begin{equation}\label{eq:2.28}
 \mathcal R=(I-\tau Q)^{-1}\ge O,\qquad
 \mathcal K=\mathcal R(I+\tau \mathcal C)\ge O,\qquad 
 \norm{\mathcal R}_{\infty\to\infty}=\norm{\mathcal K}_{\infty\to\infty}=1.
\end{equation}
\end{lemma}
\begin{proof}
By the same argument as in the proof of Theorem~\ref{Theorem:2.1},
we have $\mathcal R\ge O$ and $\mathcal R\one=\one$.
To prove positivity of the error propagation matrix~$\mathcal K$, put
$\theta=(1+\delta_Q)^{-1}$ and
$P_\theta=(1-\theta)I+\tau(\mathcal C+\theta Q)$.
For $i\ne j$ and $q_{ij}>0$, the definition of $\zeta$ and $\zeta\le\theta$ give
$(P_\theta)_{ij}=\tau(\mathcal C_{ij}+\theta q_{ij})\ge\tau(\theta-\zeta)q_{ij}\ge0$,
while $q_{ij}=0$ implies $\mathcal C_{ij}=0$ and hence $(P_\theta)_{ij}=0$.
The diagonal entries of $P_\theta$ satisfy
\begin{equation*}
 1-\theta\left(1+\tau\sum_{j\ne i}q_{ij}\right)
 \ge 1-\theta(1+\delta_Q)=0.
\end{equation*}
Hence $P_\theta\ge O$ and
$I+\tau \mathcal C=\theta(I-\tau Q)+P_\theta$, giving
$\mathcal K=\theta I+\mathcal R P_\theta\ge O$.
Moreover, $\mathcal K\one=\mathcal R(I+\tau \mathcal C)\one=\one$. Thus both $\mathcal R$ and $\mathcal K$ are nonnegative matrices with unit row sums, and $\norm{\mathcal R}_{\infty\to\infty}=\norm{\mathcal K}_{\infty\to\infty}=1$.
\qedhere
\end{proof}

To express the hypotheses of Lemma~\ref{Lemma:2.7} in terms of the mesh, define
\begin{equation}\label{eq:2.29}
 \Lambda_h=\max_{i\in\mathcal I_h,j\ne i,a_{ij}>0}
               \frac{|\ell_{ij}|}{a_{ij}},\qquad
 D_h^A=\max_i\sum_{j\ne i}a_{ij},\qquad
 \Gamma_h=\max_{i\in\mathcal I_h,K\ni x_i}\norm{E_{i,K}^{-1}}_{2\to1}.
\end{equation}
Here an empty maximum is understood to be zero. By convention, we set
$\Lambda_h=\infty$ if $a_{ij}=0$ but $\ell_{ij}\ne\mathbf0$
for some $i\in\mathcal I_h$ and $j\ne i$.
The following lemma first gives a computable contractivity condition,
then reduces it to conditions on the mesh size, time step, and gradient
bounds under~\eqref{eq:2.18} and~\eqref{eq:2.22}.
\begin{lemma}[Maximum-norm contractivity of error propagation]\label{Lemma:2.8}
Assume the sign condition~\eqref{eq:2.15} and finite
$\Lambda_h,D_h^A,\Gamma_h$ from \eqref{eq:2.29}.
Let $W,Y,Z$ satisfy \eqref{eq:2.25}, and let
$\mathcal C=\mathcal C(W,Y;Z)$ be the matrix from Lemma~$\mathrm{\ref{Lemma:2.6}}$.
 \begin{enumerate}[label=\textup{(\roman*)}]
\item For $\tau>0$, the matrices $\mathcal R=(I-\tau Q(W))^{-1}$ and
$\mathcal K=\mathcal R(I+\tau \mathcal C)$ satisfy \eqref{eq:2.28} whenever
\begin{equation*}
 \mathcal S_h(\tau;M_1,M_2):=(p-1)M_1^{p-2}M_2\Lambda_h
       [1+\tau(D_h^A+\Gamma_h M_1^{p-1})]\le1.
\end{equation*}
 \item  If, in addition, the mesh condition \eqref{eq:2.18} and the condition
\eqref{eq:2.22} hold, then, for $0\le M_2\le M_1$,
\begin{equation*}
 \mathcal S_h(\tau;M_1,M_2)
 \le C ( hM_1^{p-1}+\tau h^{-1}M_1^{p-1}+\tau M_1^{2p-2} ),
\end{equation*}
with $C$ independent of $h$, $\tau$, $M_1$, and $M_2$.
 \item Under the additional mesh assumptions in \textup{(ii)}, for
$0\le M_2\le M_1$ and $\dt>0$, which may depend on $h$,
\begin{equation}\label{eq:2.30}
 hM_1^{p-1}+\dt h^{-1}M_1^{p-1}\to0
\end{equation}
as $h\downarrow0$ implies $\mathcal S_h(\tau;M_1,M_2)\le1$ for every $0<\tau\le\dt$ and all sufficiently
small $h$.
 \end{enumerate}
 \end{lemma}
\begin{proof}
Set $Q=Q(W)$ and define $\delta_Q,\zeta$ by \eqref{eq:2.27}.
By \eqref{eq:2.9}, \eqref{eq:2.12}, and~\eqref{eq:2.15}, $Q\one=\mathbf0$ and
$q_{ij}\ge a_{ij}\ge0$ for $i\ne j$.
If $i\ne j$ and $q_{ij}(W)=0$, then $a_{ij}=0$.
Since $\Lambda_h<\infty$, this implies $\ell_{ij}=\mathbf0$.
Lemma~\ref{Lemma:2.6} therefore gives $\mathcal C_{ij}=0$ on these entries and
$\zeta\le(p-1)M_1^{p-2}M_2\Lambda_h$.
By \eqref{eq:2.6} and \eqref{eq:2.7} with $v=0$, we have
$|b_i(W)|=|g_i(W)|^{p-1}\le M_1^{p-1}$.
Thus \eqref{eq:2.9}, \eqref{eq:2.10}, \eqref{eq:2.12}, and the definitions of $\delta_Q$, $D_h^A$, and $\Gamma_h$ in~\eqref{eq:2.27} and~\eqref{eq:2.29} give
$\delta_Q\le\tau(D_h^A+\Gamma_h M_1^{p-1})$.
Consequently, $\mathcal S_h(\tau;M_1,M_2)\le1$ implies $\zeta(1+\delta_Q)\le1$, and
Lemma~\ref{Lemma:2.7} yields \eqref{eq:2.28}.
Under \eqref{eq:2.18} and \eqref{eq:2.22}, we have
$\Lambda_h\le c_Dh$, while Proposition~\ref{Proposition:2.3} gives
$D_h^A\lesssim h^{-2}$. Since $|z|_1\le\sqrt N|z|$, we have
$\norm{X}_{2\to1}\le\sqrt N\norm{X}_{2\to2}$, so \eqref{eq:2.19}
implies $\Gamma_h\lesssim h^{-1}$.
Substitution into the definition of $\mathcal S_h$, with $M_2\le M_1$,
proves the bound in \textup{(ii)}. Its last term is the product of the first two,
so \eqref{eq:2.30} makes the bound tend to zero uniformly for
$0<\tau\le\dt$.
\qedhere
\end{proof}


\subsection{Interpolation, inverse, and Ritz estimates}\label{Subsection:2.3}
The error analysis uses inverse estimates to bound numerical gradients
from nodal errors and a Ritz projection to control the diffusion residual.
We first record the Poincar\'e, inverse, and interpolation estimates.
We then construct the projection on a fixed
outer domain, so that its estimates remain uniform as $\Omega_h$ varies.
\begin{lemma}[Poincar\'e inequalities]\label{Lemma:2.9}
Let $E\subset\R^N$ be a fixed bounded Lipschitz domain.
There is a constant $C>0$, depending only on $N$ and an upper bound for
$|E|$, such that
\begin{gather}
 \norm{z}_{2,E}\le\norm{z}_{H^1(E)}
 \le C\norm{\nabla z}_{2,E},
 \qquad
 z\in H_0^1(E),\label{eq:2.31}\\
 \norm{z}_{\infty,E}\le\norm{z}_{W^{1,\infty}(E)}
 \le C\norm{\nabla z}_{\infty,E},
 \qquad
 z\in W^{1,\infty}(E)\cap H_0^1(E).\label{eq:2.32}
\end{gather}
\end{lemma}
\begin{proof}
See \cite{BrennS2008}*{Proposition~5.3.5} and
\cite{GilbaT1983}*{Section~7.8}, letting $q\to\infty$ in the latter
for the $L^\infty$ estimate.
\end{proof}

\begin{lemma}[Inverse and interpolation estimates]\label{Lemma:2.10}
Let $\Ih$ be the nodal operator \eqref{eq:2.3}, and assume the mesh
condition \eqref{eq:2.18}.
For every $w_h\in V_h$, with nodal vector $W=(W_i)$, $W_i=w_h(x_i)$,
\begin{equation}\label{eq:2.33}
 \norm{\nabla w_h}_{\infty,\Omega_h}\le C_{\rm inv}h^{-1}\norm{W}_{\ell^\infty}.
\end{equation}
Moreover, for $f\in W^{1,\infty}(\Omega_h)$, using its continuous representative,
\begin{equation}\label{eq:2.34}
 \norm{\nabla\Ih f}_{\infty,\Omega_h}\le C\norm{\nabla f}_{\infty,\Omega_h}.
\end{equation}
For $f\in W^{2,\infty}(\Omega_h)$, we also have
\begin{equation}\label{eq:2.35}
 \norm{\nabla(f-\Ih f)}_{\infty,\Omega_h}\le Ch\norm{D^2f}_{\infty,\Omega_h}.
\end{equation}
All constants depend only on $N$ and the uniform mesh constants in
\eqref{eq:2.18}, and are uniform over the varying domains $\Omega_h$.
\end{lemma}
\begin{proof}
Estimate~\eqref{eq:2.33} follows from~\cite{BrennS2008}*{Theorem~4.5.11}.
The interpolation estimates~\eqref{eq:2.34} and~\eqref{eq:2.35} follow from~\cite{BrennS2008}*{Theorem~4.4.20}, applied elementwise.
\qedhere
\end{proof}

We next construct finite element approximations of the continuous solution
to serve as reference functions in the error analysis of
Section~\ref{Subsection:3.2}. On general simplicial meshes, the
mass-lumped discrete Laplacian of the nodal interpolant need not be
pointwise consistent with the continuous Laplacian. We therefore use the
classical Ritz projection, whose defining orthogonality identifies the discrete
Laplacian with weighted nodal averages of the continuous Laplacian.
This identity allows us to use the continuous equation to estimate the
residual; see~\cite{Thome2006}.
To handle the varying domains $\Omega_h$, we extend the continuous
solution to a fixed outer domain, take its Ritz projection there, and
restrict the result to $\Omega_h$. This construction adapts those in
\cites{Schat1980,JenseS2013,LeykeL2021,GaoQ2021,LiQXY2024}.
For related finite element analyses on approximate domains, see~\cites{BrambDT1972,KashiK2020,Tabat1978}.
For $L^\infty$ and $W^{1,\infty}$ error estimates on approximate domains, see~\cites{GongLX2023,KashiK2020a}.

We now specify the domains and meshes used in this construction.
Let $\Omega$ be a fixed bounded $C^2$ domain, and let $D\subset\R^N$ be a
fixed bounded convex polytope. Assume that the computational meshes are
contained in conforming simplicial meshes $\widehat{\mathcal T}_h$ of $D$,
with
\begin{equation}\label{eq:2.36}
 \overline\Omega\subset D,\qquad
 \mathcal T_h\subset\widehat{\mathcal T}_h,\qquad
 \widehat h:=\max_{K\in\widehat{\mathcal T}_h}\operatorname{diam}K\simeq h.
\end{equation}
The background meshes $\widehat{\mathcal T}_h$ are assumed to be uniformly
shape regular and quasiuniform as in \eqref{eq:2.18}. By \eqref{eq:2.36},
every computational element is also an element of the background mesh.
The sign condition~\eqref{eq:2.15} and the condition
\eqref{eq:2.22} are imposed only on $\mathcal T_h$.
To transfer functions from $\Omega$ to $D$, fix a linear
extension operator $\mathcal E$, with $\mathcal E w=w$ on $\Omega$ and support in a fixed
compact subset of $D$, such that, for $k=1,2$ and $q=2,\infty$,
\begin{equation}\label{eq:2.37}
 \norm{\mathcal E w}_{W^{k,q}(D)}\le C\norm{w}_{W^{k,q}(\Omega)}.
\end{equation}
Such an operator is given by the extension construction in
\cite{Evans1998}*{Section~5.4, Theorem~1 and Remark~(i)}.
The support condition ensures that $\mathcal E w$ has zero trace on $\partial D$.
In particular, \eqref{eq:2.37} with $k=1$ and $q=2$ shows that
$\mathcal E:H^1(\Omega)\to H_0^1(D)$ is a bounded linear operator.

Let $\widehat V_h^0$ be the continuous piecewise affine space on
$\widehat{\mathcal T}_h$ with zero trace on $\partial D$.
For $w\in H_0^1(D)$, define its Ritz projection
$\Rh^D w\in\widehat V_h^0$ by
\begin{equation}\label{eq:2.38}
 \int_D\nabla(w-\Rh^D w)\cdot\nabla\chi_h\, \mathrm{d}x=0,
 \qquad \chi_h\in\widehat V_h^0.
\end{equation}
Equation~\eqref{eq:2.38} has a unique solution and defines a linear
projection $\Rh^D$ onto $\widehat V_h^0$; see
\cite{BrennS2008}*{Proposition~5.3.5 and Corollary~2.7.13}.
We test~\eqref{eq:2.38} with $\chi_h=\Rh^D w$ and then apply the
Cauchy--Schwarz inequality to obtain
\begin{equation*}
 \norm{\nabla\Rh^D w}_{2,D}^2
 =\int_D\nabla w\cdot\nabla\Rh^D w\,\mathrm{d}x
 \le\norm{\nabla w}_{2,D}\norm{\nabla\Rh^D w}_{2,D}.
\end{equation*}
Consequently, \eqref{eq:2.31} applied on $D$ yields
\begin{equation}\label{eq:2.39}
 \norm{\Rh^D w}_{H^1(D)}\le C\norm{w}_{H^1(D)},
 \qquad w\in H_0^1(D).
\end{equation}
Thus $\Rh^D:H_0^1(D)\to H_0^1(D)$
is a bounded linear operator.
In Section~\ref{Subsection:3.2}, we use
$(\Rh^D \mathcal E u(t))|_{\Omega_h}$ as a reference function.
To estimate this reference function, we first establish two stability bounds.
\begin{lemma}[Gradient stability]\label{Lemma:2.11}
Let $N\in\{1,2,3\}$, and assume the domain, extension, and background mesh
conditions above. There are constants $C,C_E>0$, independent of $h$,
such that, for $w\in W^{1,\infty}(D)\cap H_0^1(D)$ and
$w\in W^{1,\infty}(\Omega)\cap H_0^1(\Omega)$, respectively,
\begin{equation}\label{eq:2.40}
 \norm{\nabla\Rh^D w}_{\infty,D}
 \le C\norm{\nabla w}_{\infty,D},
 \qquad
 \norm{\nabla\Rh^D \mathcal E w}_{\infty,D}
 \le C_E\norm{\nabla w}_{\infty,\Omega}.
\end{equation}
\end{lemma}
\begin{proof}
For $N=2,3$, the first inequality in \eqref{eq:2.40} follows from
\cite{DieniRS2024}*{Theorem~3.1 and Corollary~3.2}.
For $N=1$, $\Rh^D$ equals nodal interpolation; see~\cite{Thome2006}*{Chapter~5}.
Thus the first inequality in \eqref{eq:2.40} follows directly with $C=1$.
Combining \eqref{eq:2.37} and the first inequality in \eqref{eq:2.40}
with \eqref{eq:2.32}
for functions with zero trace on $\partial\Omega$ yields
the second inequality in \eqref{eq:2.40}.
\end{proof}

We next collect the Ritz approximation estimates needed in
Section~\ref{Subsection:3.2}.

\begin{lemma}[Ritz approximation estimates]\label{Lemma:2.12}
Assume the hypotheses of Lemma~$\mathrm{\ref{Lemma:2.11}}$.
For every $w\in W^{2,\infty}(D)\cap H_0^1(D)$,
\begin{equation}\label{eq:2.41}
 \norm{\Rh^D w-w}_{W^{1,\infty}(D)}
 \le Ch\norm{D^2w}_{\infty,D}.
\end{equation}
Moreover, if $N=2$, then for every
$w\in W^{1,\infty}(D)\cap H^2(D)\cap H_0^1(D)$,
\begin{equation}\label{eq:2.42}
 \norm{w-\Rh^D w}_{\infty,D}
 \le Ch\norm{\nabla w}_{\infty,D}^{1/2}
         \norm{w}_{H^2(D)}^{1/2}.
\end{equation}
\end{lemma}
\begin{proof}
For $N=2,3$, the best approximation estimate of
\cite{DemloLSW2012}*{Theorem~2} applies on $D$, since the
quasiuniform background meshes satisfy Mesh Condition~1 of that reference.
For $N=1$, $\Rh^D$ equals nodal interpolation, as noted in the proof of
Lemma~\ref{Lemma:2.11}.
In both cases, estimate~\eqref{eq:2.35}, applied to the background mesh on $D$ using~\eqref{eq:2.36}, and~\eqref{eq:2.32} for functions with zero trace on $\partial D$ give~\eqref{eq:2.41}.
For~$N=2$, the~$L^2$ error estimate in~\cite{RannaS1982}*{Corollary}, with~$p=2$, gives
\begin{equation*}
 \norm{w-\Rh^D w}_{2,D}\le Ch^2\norm{w}_{H^2(D)}.
\end{equation*}
Applying the interpolation inequality in~\cite{Niren1959}*{Lecture~II, Theorem} to the zero extension of~$w-\Rh^D w$ to~$\R^2$, and using this estimate together with~\eqref{eq:2.40}, we obtain
\begin{equation*}
 \norm{w-\Rh^D w}_{\infty,D}
 \le C\norm{\nabla (w-\Rh^D w)}_{\infty,D}^{1/2}\norm{w-\Rh^D w}_{2,D}^{1/2}
 \le Ch\norm{\nabla w}_{\infty,D}^{1/2}\norm{w}_{H^2(D)}^{1/2},
\end{equation*}
which proves \eqref{eq:2.42}.
\end{proof}

Finally, we relate the fixed-domain projection to the discrete Laplacian
on $\Omega_h$. Define the weighted nodal averages by
\begin{equation*}
 (\Pi_h f)_i=\frac{1}{m_i}\int_{\Omega_h}f\varphi_i\, \mathrm{d}x,
 \qquad i\in\mathcal I_h.
\end{equation*}
\begin{lemma}[Nodal averages and discrete Laplacian]\label{Lemma:2.13}
Assume the hypotheses of Lemma~$\mathrm{\ref{Lemma:2.11}}$.
For every $i\in\mathcal I_h$,
\begin{gather}
 |(\Pi_h f)_i|\le\norm{f}_{\infty,\Omega_h},
 \qquad f\in L^\infty(\Omega_h),\label{eq:2.43}\\
 |(\Pi_h f)_i-f(x_i)|\le h\norm{\nabla f}_{\infty,\Omega_h},
 \qquad f\in W^{1,\infty}(\Omega_h).\label{eq:2.44}
\end{gather}
The Ritz projection also satisfies
\begin{equation}\label{eq:2.45}
 (\Delta_h((\Rh^D \mathcal E w)|_{\Omega_h}))_i=(\Pi_h\Delta w)_i,
 \qquad i\in\mathcal I_h,\quad w\in H^2(\Omega).
\end{equation}
\end{lemma}
\begin{proof}
Estimate~\eqref{eq:2.43} follows from
the nonnegativity of $\varphi_i$ and the definition of $m_i$ in~\eqref{eq:2.4}.
Next, for $f\in W^{1,\infty}(\Omega_h)$, each element in the support of
$\varphi_i$ contains $x_i$ and has diameter at most $h$. Integrating
$|f(x)-f(x_i)|\le h\norm{\nabla f}_{\infty,\Omega_h}$ against
$\varphi_i/m_i$ gives \eqref{eq:2.44}.
Finally, to prove \eqref{eq:2.45}, fix $i\in\mathcal I_h$.
Every element in the patch at $x_i$ belongs to $\mathcal T_h$.
The extension of $\varphi_i$ by zero to $D$ is therefore
the corresponding background nodal basis function. Testing
\eqref{eq:2.38} with this function and integrating by parts, we obtain
\eqref{eq:2.45} by the definition \eqref{eq:2.5} of $\Delta_h$.
\end{proof}

\section{Convergence before gradient blow-up}\label{Section:3}

On polytopal inner approximations of a smooth domain, we prove maximum-norm convergence and a uniform bound for the numerical gradients on each fixed interval before blow-up. We first combine contractivity with residual estimates for the Ritz reference function to obtain a conditional error estimate. We then use this estimate to bound the numerical gradients.

\subsection{Reference-function estimates}\label{Subsection:3.1}
We first compare $U_h^n$ with a general time-dependent finite element
function $v_h$, where $v_h(\cdot,t)\in V_h$. 
Write $V_i(t)=v_h(x_i,t)$ for its value at node $x_i$ and
$V(t)=(V_i(t))_{i\in\mathcal J_h}$ for the vector of these values,
including boundary nodes. We assume that $V(t)$ is continuous in time and allow
the boundary values $V_i(t)$, $i\in\mathcal B_h$, to be nonzero.
Define the residual on all nodes, with $Q$ as in \eqref{eq:2.11}, by
\begin{equation}\label{eq:3.1}
 r^{n+1}=\frac{V(t_{n+1})-V(t_n)}{\tau_n}
                         -Q(V(t_n))V(t_{n+1}).
\end{equation}
Since $Q$ has zero boundary rows, \eqref{eq:3.1} gives
\begin{equation}\label{eq:3.2}
 r_i^{n+1}=\frac{V_i(t_{n+1})-V_i(t_n)}{\tau_n}
 \quad\mbox{for }i\in\mathcal B_h.
\end{equation}
For $0\le T<T^*$, define the reference-gradient bound, accumulated
residual and initial discrepancy, reference approximation error, and
error outside the computational domain by
\begin{equation}
\label{eq:3.3}
\begin{aligned}
 \mathcal G_h(T)&:=\sup_{0\le t\le T}
                  \norm{\nabla v_h(\cdot,t)}_{\infty,\Omega_h},\\
 \mathcal D_h(T)&:=\norm{U^0-V(0)}_{\ell^\infty}
       +\sum_{t_{k+1}\le T}\tau_k\norm{r^{k+1}}_{\ell^\infty},\\
 \varepsilon_h(T)&:=\sup_{0\le t\le T}
                  \norm{v_h(\cdot,t)-u(t)}_{\infty,\Omega_h},\\
 \beta_h(T)&:=\sup_{0\le t\le T}
                  \norm{u(t)}_{\infty,\Omega\setminus\Omega_h}.
\end{aligned}
\end{equation}
Set $\beta_h(T)=0$ if $\Omega_h=\Omega$.
The following lemma controls the error relative to $v_h$ in terms of the initial discrepancy and the residual.

\begin{lemma}[Reference-function error estimate]\label{Lemma:3.1}
Assume the sign condition~\eqref{eq:2.15}, the mesh condition~\eqref{eq:2.18}, and finite
$\Lambda_h,D_h^A,\Gamma_h$ from \eqref{eq:2.29}.
Let $U^n$ be generated by exact solves of \eqref{eq:2.14}, with zero
numerical boundary values, and let $V(t)$ be the continuous reference
vector defined above. Fix $T<T^*$, $M_*>0$, and $M_2\ge0$ such that
$\mathcal G_h(T)\le M_2$ and $\mathcal S_h(\dt;M_*+M_2,M_2)\le1$.
If $t_n\le T$ and $G_h^k<M_*$ for $0\le k<n$, then
\begin{equation}\label{eq:3.4}
 \norm{U^n-V(t_n)}_{\ell^\infty}\le\mathcal D_h(t_n).
\end{equation}
Consequently,
\begin{equation}\label{eq:3.5}
 \norm{U_h^n-u(t_n)}_{\infty,\Omega_h}
  \le\mathcal D_h(T)+\varepsilon_h(T),\qquad
 G_h^n\le \mathcal G_h(T)+C_{\rm inv}h^{-1}\mathcal D_h(T).
\end{equation}
\end{lemma}
\begin{proof}
Set $e^k=U^k-V(t_k)$. For $k<n$, put $W=U^k$, $Y=V(t_k)$,
and $Z=V(t_{k+1})$.
Let $\mathcal C=\mathcal C(W,Y;Z)$ be the matrix
from Lemma~\ref{Lemma:2.6}.
Then, subtracting \eqref{eq:3.1} from
\eqref{eq:2.14} and applying that lemma gives
\begin{equation}\label{eq:3.6}
 (I-\tau_kQ(W))e^{k+1}=(I+\tau_k\mathcal C)e^k-\tau_kr^{k+1},
 \qquad (Q(W)-Q(Y))Z=\mathcal Ce^k.
\end{equation}
The assumptions allow us to apply Lemma~\ref{Lemma:2.8} with
$M_1=M_*+M_2$, giving \eqref{eq:2.28}. Applying the resulting
maximum-norm bounds to \eqref{eq:3.6}, we obtain
\begin{equation*}
 \norm{e^{k+1}}_{\ell^\infty}
 \le\norm{e^k}_{\ell^\infty}
       +\tau_k\norm{r^{k+1}}_{\ell^\infty},\qquad 0\le k<n.
\end{equation*}
Summing this inequality over $0\le k<n$ proves \eqref{eq:3.4}.
Finally, \eqref{eq:3.3}--\eqref{eq:3.4} give
the first bound in \eqref{eq:3.5}. For the second, the inverse inequality
\eqref{eq:2.33} gives
\begin{align*}
 G_h^n&\le \norm{\nabla v_h(\cdot,t_n)}_{\infty,\Omega_h}
       +\norm{\nabla v_h(\cdot,t_n)-\nabla U_h^n}_{\infty,\Omega_h}\\
 &\le \mathcal G_h(T)+C_{\rm inv}h^{-1}\norm{e^n}_{\ell^\infty}
 \le \mathcal G_h(T)+C_{\rm inv}h^{-1}\mathcal D_h(T).
 \qedhere
\end{align*}
\end{proof}

\subsection{Maximum-norm error estimates}\label{Subsection:3.2}

We now choose the Ritz projection from Section~\ref{Subsection:2.3} as
the reference function. To use the error bound on intervals approaching
blow-up, we track the growth of its residual in terms of the maximum
gradient norm of the exact solution up to time $T$.
For $T\in[0,T^*)$, define
\begin{equation}\label{eq:3.7}
 M(T):=\max\left\{1,\max_{0\le t\le T}\norm{\nabla u(t)}_{\infty,\Omega}\right\}.
\end{equation}
To estimate the residual $r^{n+1}$, we impose the following regularity assumption.

\par\medskip
{\setlength{\fboxsep}{8pt}
\noindent\fbox{%
\begin{minipage}{\dimexpr\linewidth-2\fboxsep-2\fboxrule\relax}
\hypertarget{regularity-R}{\textbf{Regularity assumption (R).}}
We assume that $\Omega$ is a bounded $C^2$ domain,
$u_0\in W^{2,\infty}(\Omega)\cap H_0^1(\Omega)$,
$\partial_tu$ has zero boundary trace, and
\begin{equation}\label{eq:3.8}
 \partial_t^2 u-\Delta\partial_t u
 =p|\nabla u|^{p-2}\nabla u\cdot\nabla\partial_t u.
\end{equation}
Here \eqref{eq:3.8} is understood in the sense of distributions. We also assume that
there is a constant $C_{\rm reg}>0$, depending on $\Omega,p,u_0$ but
independent of $T<T^*$, such that, with $M=M(T)$,
\begin{equation}\label{eq:3.9}
 \sup_{0\le t\le T}\left(
 \frac{\norm{D^2u(t)}_{\infty,\Omega}}{M^p}
 +\frac{\norm{\nabla \partial_t u(t)}_{\infty,\Omega}}{M^{p-1}}
 +\frac{\norm{\partial_t^2 u(t)}_{\infty,\Omega}}{M^{2p-2}}
 \right)\le C_{\rm reg}.
\end{equation}
For $N=3$, we also require
\begin{equation}\label{eq:3.10}
 \sup_{0\le t\le T}
 \frac{\norm{D^2\partial_t u(t)}_{\infty,\Omega}}{M(T)^{2p-1}}
       \le C_{\rm reg}.
\end{equation}
\end{minipage}%
}\par}
\bigskip

The following proposition gives sufficient conditions for
assumption~\hyperlink{regularity-R}{(R)}.

\begin{proposition}\label{Proposition:3.2}
Let $\Omega\subset\R^N$, $N\ge1$, be a smooth bounded domain,
let $p>2$, and let $u_0\in C_0^5(\Omega)$. Then the classical
solution of \eqref{eq:1.1} satisfies assumption~\hyperlink{regularity-R}{$\mathrm{(R)}$}.
\end{proposition}

The bound for $\nabla\partial_tu$ in \eqref{eq:3.9} near finite-time
gradient blow-up on smooth bounded domains is already available in
Porretta and Souplet
\cite{PorreS2020}*{Section~4, proof of Theorem~4.1}.
In Appendix~\ref{Subsection:A.1}, we prove all the estimates required
for the proposition, including this bound.

Under assumption~\hyperlink{regularity-R}{(R)}, \eqref{eq:1.1} and
\eqref{eq:3.9} imply
\begin{equation}\label{eq:3.11}
 \sup_{0\le t\le T}\norm{\partial_t u(t)}_{\infty,\Omega}
 \le\sqrt N\sup_{0\le t\le T}\norm{D^2u(t)}_{\infty,\Omega}+M(T)^p
 \le CM(T)^p.
\end{equation}
Together with \eqref{eq:3.9}, this bounds the time derivatives
of $u$, $\nabla u$, and $\partial_tu$ in $L^\infty(\Omega\times(0,T))$.
The fundamental theorem of calculus for $L^2(\Omega)$-valued functions
\cite{Evans1998}*{Section~5.9.2, Theorem~2} then gives the integral identities
\begin{align}
 u(t)-u(s)&=\int_s^t \partial_t u(r)\,\mathrm{d}r,\label{eq:3.12}\\
 \nabla u(t)-\nabla u(s)&=\int_s^t\nabla \partial_t u(r)\,\mathrm{d}r,\label{eq:3.13}\\
 \partial_t u(t)-\partial_t u(s)&=\int_s^t \partial_t^2 u(r)\,\mathrm{d}r,\label{eq:3.14}
\end{align}
for $0\le s\le t\le T$. These identities and the bounds in~\eqref{eq:3.9}
and~\eqref{eq:3.11} show that the resulting representatives are Lipschitz
continuous in time with values in $L^\infty(\Omega)$. We use these
representatives below.
For later use, we also derive an $H^2(\Omega)$ bound for $\partial_tu$
when $N=2$. In this case, the differentiated equation \eqref{eq:3.8} and
\eqref{eq:3.9} give
$\norm{\Delta\partial_tu(t)}_{\infty,\Omega}\le CM(T)^{2p-2}$.
The $H^2$ estimate for the homogeneous Dirichlet Poisson problem on the
fixed $C^2$ domain $\Omega$
\cite{GilbaT1983}*{Theorem~9.15 and Lemma~9.17} therefore yields
\begin{equation}\label{eq:3.15}
 \sup_{0\le t\le T}\norm{\partial_tu(t)}_{H^2(\Omega)}
 \le CM(T)^{2p-2},\qquad N=2.
\end{equation}
These bounds allow us to estimate the Ritz reference function and its
residual, including the residual at computational boundary nodes where
the reference function need not vanish.

\begin{lemma}
\label{Lemma:3.3}
Assume the hypotheses of Lemma~$\mathrm{\ref{Lemma:2.11}}$.
Suppose $T^*<\infty$,
regularity assumption~\hyperlink{regularity-R}{$\mathrm{(R)}$} holds, and
\eqref{eq:2.1} is satisfied. Set
\begin{equation}\label{eq:3.16}
 v_h(\cdot,t)=(\Rh^D \mathcal E u(t))|_{\Omega_h},\qquad V_i(t)=v_h(x_i,t) 
 \quad\mbox{for } i\in\mathcal J_h,
\end{equation}
and initialize $U_h^0$ as in Section~$\mathrm{\ref{Section:2}}$.
Then
\begin{equation}\label{eq:3.17}
 \norm{\nabla U_h^0}_{\infty,\Omega_h}
 \le C\norm{\nabla u_0}_{\infty,\Omega}.
\end{equation}
With $M=M(T)$, we also have
\begin{equation}\label{eq:3.18}
 \norm{r^{n+1}}_{\ell^\infty}
 \le C\bigl(hM^{2p-1}+\tau_nM^{2p-2}\bigr),
 \qquad t_{n+1}\le T<T^*.
\end{equation}
Moreover, for $0<h\le1$,
\begin{equation}\label{eq:3.19}
 \mathcal G_h(T)\le C_EM,\ \mathcal D_h(T)\le C(hM^{2p-1}+\dt M^{2p-2}),\ \varepsilon_h(T)\le ChM^p,\  \beta_h(T)\le ChM.
\end{equation}
\end{lemma}
\begin{proof}
The extension bound
\eqref{eq:2.37}, the second inequality in
\eqref{eq:2.40}, the approximation estimate \eqref{eq:2.41}, and
assumption~\hyperlink{regularity-R}{(R)} give
\begin{equation}\label{eq:3.20}
 \norm{\nabla v_h(\cdot,t)}_{\infty,\Omega_h}\le C_EM,\quad 
 \norm{v_h(\cdot,t)-u(t)}_{W^{1,\infty}(\Omega_h)}
 \le Ch\norm{\mathcal E u(t)}_{W^{2,\infty}(D)}\le ChM^p.
\end{equation}
Here $\norm{u(t)}_{\infty,\Omega}\le CM$ follows from
\eqref{eq:2.32} and the zero trace on $\partial\Omega$.
In particular, estimates~\eqref{eq:3.20} give the bounds
for $\mathcal G_h$ and $\varepsilon_h$.

To estimate the residual~$r^{n+1}$ for $t_{n+1}\le T$, we introduce the time average
\begin{equation*}
 w_n=\frac{1}{\tau_n}\int_{t_n}^{t_{n+1}}\partial_tu(s)\,\mathrm{d}s.
\end{equation*}
The zero trace of $\partial_tu$ and \eqref{eq:3.9} give
$w_n=0$ on $\partial\Omega$ and
$\norm{\nabla w_n}_{\infty,\Omega}\le CM^{p-1}$.
Thus \eqref{eq:2.32} and \eqref{eq:2.37} yield
\begin{equation}\label{eq:3.21}
 \norm{\nabla\mathcal E w_n}_{\infty,D}
 \le C\norm{w_n}_{W^{1,\infty}(\Omega)}
 \le C\norm{\nabla w_n}_{\infty,\Omega}\le CM^{p-1}.
\end{equation}
To express the time difference of $v_h$ in terms of $w_n$, we use
\eqref{eq:3.13} to interpret the integral identity~\eqref{eq:3.12}
in $H^1(\Omega)$. By \eqref{eq:2.37} with $k=1$, $q=2$ and
\eqref{eq:2.39}, the operators
$\mathcal E:H^1(\Omega)\to H_0^1(D)$ and
$\Rh^D:H_0^1(D)\to H_0^1(D)$ are bounded and linear.
Since they are independent of time, applying them to this identity
and using the definition \eqref{eq:3.16} gives
\begin{equation}\label{eq:3.22}
 \frac{v_h(\cdot,t_{n+1})-v_h(\cdot,t_n)}{\tau_n}
       =(\Rh^D \mathcal E w_n)|_{\Omega_h}.
\end{equation}
To compare the right-hand side of \eqref{eq:3.22} with the time average
$w_n$, we first prove the following bound at all computational nodes:
\begin{equation}\label{eq:3.23}
 \max_{i\in\mathcal J_h}|(\Rh^D \mathcal E w_n)(x_i)-w_n(x_i)|
       \le ChM^{2p-1}.
\end{equation}
For $N=1$, this difference is zero because $\Rh^D$ equals background
nodal interpolation.
For $N=2$, \eqref{eq:3.15} gives
$\norm{w_n}_{H^2(\Omega)}\le CM^{2p-2}$.
The extension estimate \eqref{eq:2.37} therefore yields
\(
 \norm{\mathcal E w_n}_{H^2(D)}
 \le C\norm{w_n}_{H^2(\Omega)}\le CM^{2p-2}.
\)
Since $\mathcal E w_n\in W^{1,\infty}(D)\cap H^2(D)\cap H_0^1(D)$
and $\mathcal E w_n=w_n$ on $\Omega$, combining this bound with
\eqref{eq:2.42} and \eqref{eq:3.21} gives
\begin{equation*}
 \max_{i\in\mathcal J_h}|(\Rh^D\mathcal E w_n)(x_i)-w_n(x_i)|
 \le\norm{\Rh^D\mathcal E w_n-\mathcal E w_n}_{\infty,D}
 \le ChM^\frac{p-1}{2}M^{p-1}
 \le ChM^{2p-1},
\end{equation*}
where the last inequality uses $M\ge1$ and $3(p-1)/2\le2p-1$.
For $N=3$, \eqref{eq:3.10} gives
$\norm{D^2w_n}_{\infty,\Omega}\le CM^{2p-1}$.
Together with \eqref{eq:3.21} and $M\ge1$, this implies
\begin{equation*}
 \norm{w_n}_{W^{2,\infty}(\Omega)}
 \le C\bigl(\norm{w_n}_{W^{1,\infty}(\Omega)}
             +\norm{D^2w_n}_{\infty,\Omega}\bigr)
 \le C\bigl(M^{p-1}+M^{2p-1}\bigr)\le CM^{2p-1}.
\end{equation*}
Applying \eqref{eq:2.41} and then \eqref{eq:2.37}, we obtain
\begin{equation*}
 \max_{i\in\mathcal J_h}|(\Rh^D\mathcal E w_n)(x_i)-w_n(x_i)|
 \le\norm{\Rh^D\mathcal E w_n-\mathcal E w_n}_{\infty,D}
 \le Ch\norm{\mathcal E w_n}_{W^{2,\infty}(D)}
 \le ChM^{2p-1},
\end{equation*}
which proves \eqref{eq:3.23} also for $N=3$.

With \eqref{eq:3.23} established, we now decompose the residual at
an interior node. Using \eqref{eq:3.22} to rewrite the time difference,
together with \eqref{eq:1.1}, \eqref{eq:2.6}, \eqref{eq:2.10}, \eqref{eq:2.11}, \eqref{eq:2.45},
\eqref{eq:3.1}, and \eqref{eq:3.16}, gives
\begin{equation*}
 \begin{aligned}
 r_i^{n+1}
 &= (\Rh^D \mathcal E w_n)(x_i)
       -(\Pi_h\Delta u(t_{n+1}))_i
       -\bigl(B_h(V(t_n))V(t_{n+1})\bigr)_i\\
 &= (\Rh^D \mathcal E w_n)(x_i)-(\Pi_hw_n)_i
       +(\Pi_h(w_n-\partial_tu(t_{n+1})))_i\\
 &\qquad +(\Pi_h|\nabla u(t_{n+1})|^p)_i
       -b(g_i(V(t_n)))\cdot(\nabla v_h(\cdot,t_{n+1}))_{K_i},
 \end{aligned}
\end{equation*}
where $K_i$ is the upwind element at $x_i$.
The first difference is a spatial projection error, the second is a
time-discretization error, and the last difference compares the exact
nonlinearity with the upwind drift. We estimate them in this order.
First, estimates~\eqref{eq:2.44}, \eqref{eq:3.21}, and
\eqref{eq:3.23} bound the first term,
$(\Rh^D \mathcal E w_n)(x_i)-(\Pi_hw_n)_i$, by $ChM^{2p-1}$.
Next, the integral identity \eqref{eq:3.14}, together with \eqref{eq:2.43} and \eqref{eq:3.9}, bounds the second term by $C\tau_nM^{2p-2}$.
To estimate the remaining nonlinear terms, put
\begin{equation*}
 \xi_1=\nabla u(x_i,t_{n+1}),\qquad \xi_2=g_i(V(t_n)),\qquad
 \xi_3=(\nabla v_h(\cdot,t_{n+1}))_{K_i}.
\end{equation*}
If $\xi_2=\mathbf0$, choose any element containing $x_i$ as $K_i$.
Estimate~\eqref{eq:2.7}, together with
\eqref{eq:3.9}, \eqref{eq:3.13}, and \eqref{eq:3.20}, gives
$|\xi_1|+|\xi_2|+|\xi_3|\le CM$ and
\begin{equation*}
 |\xi_3-\xi_1|\le ChM^p,\qquad
 |\xi_2-\xi_1|\le ChM^p+C\tau_nM^{p-1}.
\end{equation*}
Together with \eqref{eq:2.26}, these bounds give
\begin{equation*}
 |b(\xi_2)\cdot \xi_3-|\xi_1|^p|
 \le |b(\xi_2)|\,|\xi_3-\xi_1|+|b(\xi_2)-b(\xi_1)|\,|\xi_1|
 \le C(hM^{2p-1}+\tau_nM^{2p-2}).
\end{equation*}
Finally, \eqref{eq:3.9} gives
$\norm{\nabla(|\nabla u(t_{n+1})|^p)}_{\infty,\Omega}
\le pM^{p-1}\norm{D^2u(t_{n+1})}_{\infty,\Omega}\le CM^{2p-1}$,
so \eqref{eq:2.44} gives
\begin{equation*}
 \bigl|(\Pi_h|\nabla u(t_{n+1})|^p)_i
       -|\nabla u(x_i,t_{n+1})|^p\bigr|\le ChM^{2p-1}.
\end{equation*}
Combining these estimates in the decomposition of $r_i^{n+1}$
proves \eqref{eq:3.18} at interior nodes.
It remains to bound the residual at the computational boundary.
At each boundary node, \eqref{eq:3.2}
and \eqref{eq:3.22} give
$r_i^{n+1}=(\Rh^D \mathcal E w_n)(x_i)$. Since $w_n=0$ on $\partial\Omega$,
integration along the segment to a nearest point of $\partial\Omega$,
together with \eqref{eq:2.1} and \eqref{eq:3.21}, gives
\begin{equation*}
 |w_n(x_i)|\le Ch\norm{\nabla w_n}_{\infty,\Omega}
       \le ChM^{p-1}.
\end{equation*}
Combining this bound with \eqref{eq:3.23} proves
\eqref{eq:3.18} at boundary nodes.

Having proved \eqref{eq:3.18} on all nodes, we turn to the remaining
bounds in \eqref{eq:3.19} and the initialization bound \eqref{eq:3.17}.
First, the same boundary-segment argument as above gives $\beta_h(T)\le ChM$.
Next, to estimate the initial discrepancy, let
$Z^0=(u_0(x_i))_{i\in\mathcal J_h}$ be the nodal vector of $\Ih u_0$.
By the initialization, $U_i^0=Z_i^0$ at interior nodes and $U_i^0=0$
at boundary nodes. Since $u_0=0$ on $\partial\Omega$, \eqref{eq:2.1} gives
\begin{equation}\label{eq:3.24}
 \norm{U^0-Z^0}_{\ell^\infty}
       \le Ch\norm{\nabla u_0}_{\infty,\Omega}.
\end{equation}
Since $\nabla U_h^0=\nabla\Ih u_0+\nabla(U_h^0-\Ih u_0)$,
estimates~\eqref{eq:2.33}, \eqref{eq:2.34}, and
\eqref{eq:3.24} give
\begin{equation*}
 \norm{\nabla U_h^0}_{\infty,\Omega_h}
 \le C\norm{\nabla u_0}_{\infty,\Omega}
       +C_{\rm inv}h^{-1}\norm{U^0-Z^0}_{\ell^\infty}
 \le C\norm{\nabla u_0}_{\infty,\Omega}.
\end{equation*}
This proves \eqref{eq:3.17}.
To estimate $\mathcal D_h(T)$, we first
combine \eqref{eq:3.24} with the bound for
$\varepsilon_h$ in \eqref{eq:3.20} to obtain
\begin{equation*}
 \norm{U^0-V(0)}_{\ell^\infty}
 \le\norm{U^0-Z^0}_{\ell^\infty}+\varepsilon_h(0)\le ChM^p.
\end{equation*}
Moreover, the steps with $t_{k+1}\le T$ satisfy
\begin{equation*}
 \sum_{t_{k+1}\le T}\tau_k\le T^*,\qquad
 \sum_{t_{k+1}\le T}\tau_k^2\le\dt T^*.
\end{equation*}
Finally, summing \eqref{eq:3.18} over the steps with $t_{k+1}\le T$ and using
these initial and time-step bounds proves the bound for
$\mathcal D_h(T)$.
 \end{proof}

The residual bounds in Lemma~\ref{Lemma:3.3} now make the estimate of
Lemma~\ref{Lemma:3.1} explicit.
\begin{lemma}\label{Lemma:3.4}
Assume the hypotheses of Lemma~$\mathrm{\ref{Lemma:3.3}}$, the sign condition~\eqref{eq:2.15}, the mesh condition~\eqref{eq:2.18}, and~\eqref{eq:2.22}.
Let $u$ be a classical solution of~\eqref{eq:1.1}, and let $v_h$ and $V(t)$ be given by~\eqref{eq:3.16}.
Let $U^n$ be generated by exact solves of~\eqref{eq:2.14} on the time grids of Section~$\mathrm{\ref{Section:2}}$, with zero numerical boundary values.
Fix $T<T^*$ and $0<h\le1$, let $M=M(T)$, $M_2=C_EM$, and choose $M_*>0$ such that
\begin{equation}\label{eq:3.25}
 \mathcal S_h(\dt;M_*+M_2,M_2)\le1.
\end{equation}
If $G_h^k<M_*$ for $0\le k<n$ and $t_n\le T$, then
\begin{gather}
 \norm{U^n-V(t_n)}_{\ell^\infty}
 \le\norm{U^0-V(0)}_{\ell^\infty}
                +\sum_{k<n}\tau_k\norm{r^{k+1}}_{\ell^\infty}
 \le C(hM^{2p-1}+\dt M^{2p-2}),
 \label{eq:3.26}\\
 \norm{\overline U_h^n-u(t_n)}_{\infty,\Omega}\le
                  C(hM^{2p-1}+\dt M^{2p-2}).\label{eq:3.27}
\end{gather}
Here $r^{k+1}$ is the reference residual on all nodes and
$\overline U_h^n$ is the zero extension of $U_h^n$ to $\Omega$.
The constant $C$ is independent of $h$, $T<T^*$, and $M_*$.
\end{lemma}
\begin{proof}
Lemma~\ref{Lemma:3.3} gives $\mathcal G_h(T)\le M_2=C_EM$.
Thus Lemma~\ref{Lemma:3.1} applies under
\eqref{eq:3.25} and gives the first inequality in \eqref{eq:3.26}.
The second follows from the bound for $\mathcal D_h(T)$ in
\eqref{eq:3.19}.
Next, the reference approximation error and the error on $\Omega\setminus\Omega_h$ give
\begin{equation*}
 \norm{\overline U_h^n-u(t_n)}_{\infty,\Omega}
 \le\norm{U^n-V(t_n)}_{\ell^\infty}
       +\varepsilon_h(T)+\beta_h(T).
\end{equation*}
Combining \eqref{eq:3.19}, \eqref{eq:3.26}, and
this bound, with $M\ge1$, proves \eqref{eq:3.27}.
Since the estimates in \eqref{eq:3.19} are independent of $M_*$, so are the
constants in \eqref{eq:3.26} and \eqref{eq:3.27}.
 \end{proof}

We now prove convergence and a uniform bound for the numerical
gradients on each fixed interval before blow-up.
\begin{theorem}\label{Theorem:3.5}
Let $\Omega\subset\R^N$, $N\in\{1,2,3\}$, be a fixed bounded $C^2$ domain, and let $p>2$.
Let $u$ be a classical solution of~\eqref{eq:1.1} with $T^*<\infty$ satisfying regularity assumption~\hyperlink{regularity-R}{$\mathrm{(R)}$}.
Let the polytopal inner approximations $\Omega_h\subset\Omega$ satisfy~\eqref{eq:2.1}, and let their conforming simplicial meshes $\mathcal T_h$ satisfy the sign condition~\eqref{eq:2.15}, the mesh condition~\eqref{eq:2.18}, and~\eqref{eq:2.22}.
Assume that these meshes extend to conforming simplicial meshes $\widehat{\mathcal T}_h$ of a fixed bounded convex polytope $D$ as in~\eqref{eq:2.36}, and that the background meshes satisfy~\eqref{eq:2.18} with $\mathcal T_h$ and $h$ replaced by $\widehat{\mathcal T}_h$ and $\widehat h$, respectively.
For each $h$, use a time grid as in Section~$\mathrm{\ref{Section:2}}$, with only finitely many nodes in every bounded time interval.
Initialize $U_i^0=u_0(x_i)$ for $i\in\mathcal I_h$ and $U_i^0=0$ for $i\in\mathcal B_h$, and generate $U^n$ by exact solves of~\eqref{eq:2.14}, with zero numerical boundary values.
Let $\overline U_h^n$ denote the zero extension of $U_h^n$ to $\Omega$.
Suppose $h\downarrow0$ and $\dt h^{-1}\to0$.
Then, for every fixed $T<T^*$,
\begin{equation}\label{eq:3.28}
 \max_{t_n\le T}\norm{\overline U_h^n-u(t_n)}_{\infty,\Omega}\le C_T(h+\dt),\qquad
 \max_{t_n\le T}G_h^n\le C_T,
\end{equation}
for all sufficiently small $h$, with $C_T$ independent of $h$.
\end{theorem}
\begin{proof}
Fix $T<T^*$ and put $M=M(T)$ and $M_2=C_EM$. Let $M_*>0$ be chosen later.
Take the Ritz reference function $v_h$ and its nodal vector $V(t)$ from~\eqref{eq:3.16}, using the extension operator and Ritz projection of Section~\ref{Subsection:2.3}.
If \eqref{eq:3.25} holds, then, for every $n$ such that $t_n\le T$ and
$G_h^k<M_*$ for $0\le k<n$, Lemma~\ref{Lemma:3.4}, the
reference-gradient bound $\mathcal G_h(T)\le C_EM$ from Lemma~\ref{Lemma:3.3},
\eqref{eq:2.33}, and \eqref{eq:3.26} give
\begin{equation*}
 \begin{aligned}
 G_h^n&\le 
 \|\nabla v_h(t_n)\|_{\infty,\Omega_h} + \|\nabla (U_h^n - v_h(t_n))\|_{\infty,\Omega_h} 
 \le 
 C_EM+C_{\rm inv}h^{-1}\norm{U^n-V(t_n)}_{\ell^\infty}\\
 &\le C_EM+C ( M^{2p-1}+\dt h^{-1} M^{2p-2} )\le C_T
 \end{aligned}
\end{equation*}
for small $h$, where $C_T>0$ is a constant independent of $h$ and $M_*$.
 Since the initial gradients are uniformly bounded by
\eqref{eq:3.17}, we can choose $M_*>\max\{C_T,\sup_h G_h^0\}$.
For this fixed~$M_*$, condition~\eqref{eq:2.30} holds with $M_1=M_*+C_EM$ and $M_2=C_EM$, since $M=M(T)$ is fixed and $h+\dt h^{-1}\to0$. Hence Lemma~\ref{Lemma:2.8}\textup{(iii)} gives~\eqref{eq:3.25} for all sufficiently small~$h$. Suppose $n_*$ is the first index with
$t_{n_*}\le T$ and $G_h^{n_*}\ge M_*$.
Lemma~\ref{Lemma:3.4} and the preceding gradient bound would
then give $G_h^{n_*}\le C_T<M_*$, a contradiction.
Thus these bounds hold throughout $[0,T]$, and
\eqref{eq:3.27} proves \eqref{eq:3.28}.
\end{proof}

\subsection{Inexact linear solves}\label{Subsection:3.3}
We now allow the linear system~\eqref{eq:2.14} to be solved approximately at each time step, with the zero boundary values still imposed exactly. For the computed vector $U^{k+1}$, let $\rho^{k+1}$ denote the residual of the linear system, defined by
\begin{equation*}
 (I-\tau_kQ(U^k))U^{k+1}=U^k+\rho^{k+1}.
\end{equation*}
The boundary rows give $\rho_i^{k+1}=0$ for $i\in\mathcal B_h$. Since $\rho^{k+1}$ is defined before division by $\tau_k$, its contribution over $n$ steps is bounded by the unweighted sum $\Sigma_n=\sum_{k<n}\norm{\rho^{k+1}}_{\ell^\infty}$. Indeed, applying the contraction estimate~\eqref{eq:2.17} at each step and then the inverse inequality~\eqref{eq:2.33} gives
\begin{equation*}
 \norm{U_h^n}_{\infty,\Omega_h}
 \le\norm{U_h^0}_{\infty,\Omega_h}+{\Sigma}_n,\qquad
 \norm{\nabla U_h^n}_{\infty,\Omega_h}
 \le C_{\rm inv}h^{-1}
       (\norm{U_h^0}_{\infty,\Omega_h}+{\Sigma}_n).
\end{equation*}
The residual $\rho^{k+1}$ also appears on the right-hand side of the error equation~\eqref{eq:3.6}. Under the remaining hypotheses of Lemma~\ref{Lemma:3.1}, we apply its contraction argument at each step and sum to obtain
\begin{equation}\label{eq:3.29}
 \norm{U^n-V(t_n)}_{\ell^\infty}
 \le\norm{U^0-V(0)}_{\ell^\infty}
       +\sum_{k<n}\tau_k\norm{r^{k+1}}_{\ell^\infty}
       +{\Sigma}_n
 \le C_T(h+\dt)+{\Sigma}_n.
\end{equation}
The sum of the first two terms on the right is $\mathcal D_h(t_n)$. For fixed $T<T^*$ and $t_n\le T$, the bound for $\mathcal D_h(T)$ in~\eqref{eq:3.19} gives the last inequality under the remaining hypotheses of Theorem~\ref{Theorem:3.5}. Here $C_T$ is independent of $h$, $M_*$, and the linear solver residuals.

To ensure ${\Sigma}_n=O_T(h)$, it suffices to impose
\begin{equation}\label{eq:3.30}
 \norm{\rho^{k+1}}_{\ell^\infty}\le\eta h\tau_k,
\end{equation}
where $\eta\ge0$ is a constant independent of $h$ and $k$.
Summing \eqref{eq:3.30} gives
\begin{equation}\label{eq:3.31}
 {\Sigma}_n\le\eta h t_n.
\end{equation}
Thus the linear solver residuals contribute at most $\eta hT$ to the
maximum-norm error bound on $[0,T]$. Since
$h^{-1}{\Sigma}_n\le\eta T$, the inverse estimate \eqref{eq:2.33}
adds at most $C_{\rm inv}\eta T$ to the gradient bound used in the proof
of Theorem~\ref{Theorem:3.5}. Increasing $M_*$ accordingly, the same
proof shows that the convergence result \eqref{eq:3.28} remains valid
for inexact solves under the remaining hypotheses, with constants that
may depend on $\eta$.
\section{Approximation of the gradient blow-up time}\label{Section:4}

We prove convergence of the gradient-threshold times to~$T^*$ on polytopal inner approximations of~$\Omega$. We first establish the existence of a convergent threshold scale and then derive sufficient conditions for convergence with prescribed scales. The analysis also covers inexact linear solves with controlled residuals.

In this section, limits and asymptotic estimates for mesh-dependent
quantities are understood as $h\downarrow0$, unless another limiting
variable is specified.
Recall the maximum norm of the computed gradient on $\Omega_h$ and the
corresponding gradient-threshold time:
\begin{equation}\label{eq:4.1}
 G_h^n=\norm{\nabla U_h^n}_{\infty,\Omega_h},\qquad
 T_h(\kappa)=\inf\{t_n\colon G_h^n\ge\kappa L_h\}.
\end{equation}
Here $L_h>0$, $L_h\to\infty$, and $\kappa>0$.
As noted in Section~\ref{Section:1}, the gradient of the numerical solution
is bounded on every fixed mesh by Theorem~\ref{Theorem:2.1} and
\eqref{eq:2.33}, so diverging gradient thresholds can be reached only under
refinement. A related obstruction for conservative finite volume schemes
is noted in \cite{ZhouS2017}.


\subsection{Existence of a convergent threshold scale}

The zero extension allows us to compare solution values at fixed points as
the computational domain varies. To obtain the following estimate, we subdivide each line segment at mesh faces and apply the fundamental theorem of calculus on each subsegment.

\begin{lemma}\label{Lemma:4.1}
Let $v_h$ be a continuous $P_1$ function on a polytopal domain $\Omega_h$
triangulated by a finite simplicial mesh, with zero trace on every component
of its boundary. Define its zero extension by $\overline v_h=v_h$ on
$\Omega_h$ and $\overline v_h=0$ elsewhere in $\R^N$. Then
\begin{equation}\label{eq:4.2}
 |\overline v_h(x)-\overline v_h(y)|
 \le \norm{\nabla v_h}_{\infty,\Omega_h} |x-y|
\end{equation}
for $x,y\in\R^N$.
\end{lemma}

For $t\in[0,T^*)$, write $G(t):=\norm{\nabla u(t)}_{\infty,\Omega}$.
The following theorem shows that value convergence and uniform bounds for
the numerical gradients on each fixed interval before $T^*$ suffice to
choose a scale $L_h\to\infty$ for which the gradient-threshold times converge to $T^*$.

\begin{theorem}
\label{Theorem:4.2}
Let $\Omega\subset\R^N$ be a fixed bounded Lipschitz domain and let
$u$ be a solution of \eqref{eq:1.1} with $u$ and $\nabla u$ continuous on
$\overline\Omega\times[0,T]$ for every $T<T^*$. Suppose that $T^*$ is
finite and \eqref{eq:1.2} holds.
Let $\mathcal T_h$ be finite conforming simplicial meshes of polytopal
domains $\Omega_h\subset\Omega$, indexed by $h\downarrow0$.
For each~$h$, let $t_n=t_n(h)$ be the nodes of a time grid as in
Section~$\mathrm{\ref{Section:2}}$, with $\dt=\dt(h)\downarrow0$ as $h\downarrow0$.
Assume that, for some~$T_{\mathrm{end}}>T^*$ independent of~$h$, each time grid covers~$[0,T_{\mathrm{end}}]$.
Let $U_h^n$ be continuous $P_1$ functions
on $\Omega_h$ with zero trace on every boundary component, and let
$\overline U_h^n$ denote their zero extensions to $\Omega$.
Assume that for every fixed $T<T^*$,
\begin{equation}\label{eq:4.3}
 \max_{t_n\le T}\norm{\overline U_h^n-u(t_n)}_{\infty,\Omega}\to0,
 \qquad
 \max_{t_n\le T}G_h^n\le C_T,
\end{equation}
for all sufficiently small $h$, where $C_T$ is independent of $h$.
Then there exists a threshold scale $L_h\to\infty$ such that
\begin{equation}\label{eq:4.4}
 \sup_{\kappa\in[\kappa_-,\kappa_+]}
       |T_h(\kappa)-T^*|\to0
\end{equation}
for every $0<\kappa_-\le\kappa_+<\infty$.
\end{theorem}
\begin{proof}
Choose increasing times $s_j\uparrow T^*$ as $j\to\infty$, with
$G(s_j)>j^2$ for $j=1,2,\ldots$.
For fixed $j$, let $n_j(h)=\max\{n\colon t_n(h)\le s_j\}$, so that
$t_{n_j(h)}\le s_j<t_{n_j(h)+1}$.
Since $0\le s_j-t_{n_j(h)}<\tau_{n_j(h)}\le\dt$, \eqref{eq:4.3}
and time continuity imply
$\norm{\overline U_h^{n_j(h)}-u(s_j)}_{\infty,\Omega}\to0$ as $h\downarrow0$.
Therefore, passing to the limit in \eqref{eq:4.2} for fixed
$x,y\in\overline\Omega$ gives
\begin{equation*}
 |u(x,s_j)-u(y,s_j)|
 \le |x-y| \liminf_{h\downarrow0}G_h^{n_j(h)}.
\end{equation*}
Taking directional derivatives in this inequality at interior
points and using continuity
of $\nabla u$ up to the boundary gives
\begin{equation*}
 G(s_j)\le\liminf_{h\downarrow0}G_h^{n_j(h)}.
\end{equation*}
Since $G(s_j)>j^2$, this lower bound implies, for each fixed $j$,
\(
 G_h^{n_j(h)}>j^2
\)
for all sufficiently small $h$.

We now choose a strictly decreasing sequence $h_j\downarrow0$ such that
the preceding gradient lower bound holds for $0<h\le h_j$, and set
$L_h=j$ for $h_{j+1}<h\le h_j$. For $h>h_1$, choose any positive $L_h$.
Fix $\kappa>0$. For $h_{j+1}<h\le h_j$ with $j\ge\kappa$,
this bound and $L_h=j$ give
$G_h^{n_j(h)}>j^2\ge\kappa L_h$.
By \eqref{eq:4.1}, $T_h(\kappa)\le t_{n_j(h)}\le s_j<T^*$, proving
$\limsup_{h\downarrow0}T_h(\kappa)\le T^*$.
For the lower bound, fix $0<T<T^*$.
Since $L_h\to\infty$, the gradient bound in \eqref{eq:4.3} gives
$G_h^n\le C_T<\kappa L_h$ for every $t_n\le T$ and all sufficiently small $h$.
Thus $T_h(\kappa)>T$, and letting $T\uparrow T^*$ proves convergence.
Finally, for $\kappa\in[\kappa_-,\kappa_+]$, \eqref{eq:4.1} gives
$T_h(\kappa_-)\le T_h(\kappa)\le T_h(\kappa_+)$.
Convergence at the two endpoints proves \eqref{eq:4.4}
for the same scale $L_h$.
 \end{proof}

\begin{remark}
\label{Remark:4.3}

Under the hypotheses of Theorem~\ref{Theorem:3.5} and~\eqref{eq:1.2}, there exists a threshold scale~$L_h\to\infty$ for which~\eqref{eq:4.4} holds. Indeed, estimate~\eqref{eq:3.28} verifies~\eqref{eq:4.3}, so Theorem~\ref{Theorem:4.2} applies. This remains valid for inexact linear solves
under \eqref{eq:3.30}, with zero boundary values imposed exactly,
as shown in Section~\ref{Subsection:3.3}.

%

\end{remark}

\subsection{Convergence for prescribed thresholds}\label{Subsection:4.2}

We now derive sufficient conditions for convergence with prescribed threshold scales. We first prove a geometric lemma used in the proof.

\begin{lemma}\label{Lemma:4.4}
For a fixed bounded Lipschitz domain $\Omega\subset\R^N$, there are constants
$a_\Omega>0$ and $\ell_0>0$, depending
only on the domain, such that for every $x\in\overline\Omega$ and
$q\in\R^N$ one can choose a unit vector $e$ with
\begin{equation}\label{eq:4.5}
 x+se\in\overline\Omega,\qquad |q\cdot e|\ge a_\Omega|q|,
\end{equation}
for every $0\le s\le \ell_0$.
\end{lemma}
\begin{proof}
When $N=1$, write $\Omega=(x_-,x_+)$ and choose
$a_\Omega=1$ and $\ell_0=(x_+-x_-)/2$. From any $x\in[x_-,x_+]$, at least one
endpoint is at distance at least $\ell_0$. Choose $e\in\{-1,1\}$ toward
that endpoint. Then $x+se\in[x_-,x_+]$ for $0\le s\le \ell_0$ and
$|qe|=|q|$, proving \eqref{eq:4.5}.

For $N\ge2$, write $x=(x',x_N)$ and $e=(e',e_N)$, with
$x',e'\in\R^{N-1}$. In a boundary chart, $\Omega$ is locally given by
$x_N>\psi(x')$, where $\psi$ has Lipschitz constant $L_\psi$.
For $x\in\overline\Omega$ in this chart, choose a unit vector $e$ satisfying
$e_N>(L_\psi+1)|e'|$. For $s\ge0$ such that $x+se$ remains in the coordinate neighborhood,
the Lipschitz condition and $x_N\ge\psi(x')$ give
\begin{equation*}
 x_N+se_N-\psi(x'+se')\ge x_N-\psi(x')+s(e_N-L_\psi|e'|)\ge0.
\end{equation*}
These directions therefore form an inward cone. A finite cover by smaller
charts gives cones near the boundary with a common opening angle and
a common segment length $\ell_0>0$. Away from the boundary, interior balls
of a common positive radius give the same conclusion, after decreasing
$\ell_0$ if necessary. If $v$
is the unit cone axis and $v_1,\ldots,v_{N-1}$ span $v^\perp$ orthonormally,
choose a fixed angle $0<\alpha<\pi/2$ smaller than the common opening angle and
set $e_{j,\pm}=(\cos\alpha)v\pm(\sin\alpha)v_j$ for $1\le j\le N-1$.
These directions lie inside the cone and satisfy
\begin{equation*}
 \begin{aligned}
 \sum_{j=1}^{N-1}
 \bigl(|q\cdot e_{j,+}|^2+|q\cdot e_{j,-}|^2\bigr)
 &=2(N-1)\cos^2\alpha\,|q\cdot v|^2
 +2\sin^2\alpha\, \sum_{j=1}^{N-1} |q\cdot v_j|^2 
 \\
 &\ge2(N-1)\min\left\{\cos^2\alpha,\frac{\sin^2\alpha}{N-1}\right\}|q|^2.
 \end{aligned}
\end{equation*}
Thus, at least one squared projection is no smaller than the average of these
$2(N-1)$ terms. Choosing its direction
proves \eqref{eq:4.5} with
$a_\Omega=\min\{\cos\alpha,\sin\alpha/\sqrt{N-1}\}>0$.
 \end{proof}

For $M>\max\{1,G(0)\}$, let $s(M)$ be the first time the exact gradient
norm reaches $M$:
\begin{equation}\label{eq:4.6}
 s(M)=\inf\{t\in[0,T^*)\colon G(t)\ge M\}.
\end{equation}
By continuity of~$G$,~\eqref{eq:1.2}, and the definition of~$M(t)$ in~\eqref{eq:3.7}, we have
\begin{equation}\label{eq:4.7}
 0<s(M)<T^*,\qquad G(s(M))=M(s(M))=M.
\end{equation}
Moreover, $s(M)\uparrow T^*$ as $M\to\infty$.
Indeed, $s(M)$ is nondecreasing in $M$. If its limit
were less than $T^*$, choose $T<T^*$ larger than that limit.
Then $s(M)<T$ and $G(s(M))=M$ for arbitrarily large $M$, contradicting
the boundedness of the continuous function $G$ on the fixed interval $[0,T]$.

To compare the numerical solution at a nearby time level~$t_n$ with the exact solution at~$s(M)$, we must control the error caused by the time shift~$|t_n-s(M)|$.
We must also control how the exact gradient varies along the
segment used to detect a large numerical gradient.
For $\delta,r\ge0$, we define the corresponding moduli by
\begin{equation}\label{eq:4.8}
\begin{aligned}
 \omega_{\mathrm{time}}(\delta;T)
 &:=\sup\{|u(x,t)-u(x,s)|\colon x\in\overline\Omega,\ 0\le s,t\le T,\ 
                         |t-s|\le\delta\},\\
 \omega_{\mathrm{grad}}(r;T)
 &:=\sup\{|\nabla u(x,t)-\nabla u(y,t)|\colon 0\le t\le T,\ x,y\in\overline\Omega,\
                   |x-y|\le r,\ [x,y]\subset\overline\Omega\}.
\end{aligned}
\end{equation}
Here $[x,y]$ denotes the closed line segment joining~$x$ and~$y$.

The next theorem gives sufficient conditions for the gradient-threshold times to converge to~$T^*$ when the threshold scale is prescribed.

\begin{theorem}
\label{Theorem:4.5}
Let $\Omega\subset\R^N$ be a fixed bounded Lipschitz domain and let
$\Omega_h\subset\Omega$ be polytopal domains with conforming simplicial meshes.
Assume the sign condition~\eqref{eq:2.15}, the mesh condition~\eqref{eq:2.18}, and finite $\Lambda_h,D_h^A,\Gamma_h$ from \eqref{eq:2.29}.
Let $u$ be a classical solution of \eqref{eq:1.1}, with $u$ and
$\nabla u$ continuous on $\overline\Omega\times[0,T]$ for every $T<T^*$,
satisfying \eqref{eq:1.2}.
Assume that, for some~$T_{\mathrm{end}}>T^*$ independent of~$h$, each time grid specified in Section~$\mathrm{\ref{Section:2}}$ covers~$[0,T_{\mathrm{end}}]$.
Let $U^n$ be generated by exact solves of \eqref{eq:2.14}.
Assume zero initial boundary trace and uniformly bounded initial gradients:
$\sup_h\norm{\nabla U_h^0}_{\infty,\Omega_h}<\infty$.
Let $v_h(\cdot,t)\in V_h$ be a reference function as in
Section~$\mathrm{\ref{Subsection:3.1}}$, and use the moduli in \eqref{eq:4.8} and the domain constants $a_\Omega,\ell_0$ of Lemma~$\mathrm{\ref{Lemma:4.4}}$.
Fix $\kappa>0$ and let $L_h\to\infty$ as $h\downarrow0$.
Assume the following conditions.
\begin{enumerate}[label=\textup{(\roman*)}]
\item For each fixed $T<T^*$ there are finite constants $\mathcal G_T,C_T$,
independent of $h$, such that
\begin{equation}\label{eq:4.9}
 \mathcal G_h(T)\le \mathcal G_T,\qquad \mathcal D_h(T)\le C_T h.
\end{equation}
\item For a fixed $C_*>2/a_\Omega$, put
$M_h=C_*\kappa L_h$ and $s_h=s(M_h)$, with $s$ as in \eqref{eq:4.6}.
For all sufficiently small $h$,
\begin{equation}\label{eq:4.10}
 \mathcal S_h(\dt;\kappa L_h+\mathcal G_h(s_h),\mathcal G_h(s_h))\le1.
\end{equation}
\item There are lengths $r_h$ with $0<r_h\le \ell_0$ such that
\begin{gather}
 4\omega_{\mathrm{grad}}(r_h;s_h)\le a_\Omega M_h,\label{eq:4.11}\\
 \frac{\mathcal D_h(s_h)+\varepsilon_h(s_h)+\beta_h(s_h)
                    +\omega_{\mathrm{time}}(\dt;s_h)}{M_hr_h}
                  \to0.\label{eq:4.12}
\end{gather}
\end{enumerate}
Then $T_h(\kappa)\to T^*$, and these gradient-threshold times are finite for
sufficiently small $h$.
Moreover, if these hypotheses hold for both $\kappa=\kappa_-$ and $\kappa=\kappa_+$, where $0<\kappa_-\le\kappa_+<\infty$, with the same threshold scale~$L_h$ and numerical solution, then
\begin{equation}\label{eq:4.13}
 \sup_{\kappa\in[\kappa_-,\kappa_+]}|T_h(\kappa)-T^*|\to0.
\end{equation}
\end{theorem}
\begin{proof}
For the lower bound on $T_h(\kappa)$, fix $T<T^*$ and choose
$M_*>\max\{\sup_h\norm{\nabla U_h^0}_{\infty,\Omega_h},\mathcal G_T+C_{\rm inv}C_T\}$.
Since $s_h\to T^*$ and $\kappa L_h\to\infty$, we have $T\le s_h$ and $M_*\le\kappa L_h$ for all sufficiently small~$h$. Thus, by the monotonicity of~$\mathcal G_h$ and of~$\mathcal S_h$ in its last two arguments, condition~\eqref{eq:4.10} yields
\begin{equation*}
 \mathcal S_h(\dt;M_*+\mathcal G_h(T),\mathcal G_h(T))
 \le\mathcal S_h(\dt;\kappa L_h+\mathcal G_h(s_h),\mathcal G_h(s_h))\le1.
\end{equation*}
Now suppose that $n$ is the first index with $t_n\le T$ and $G_h^n\ge M_*$. The preceding inequality allows us to apply Lemma~\ref{Lemma:3.1} with $M_2=\mathcal G_h(T)$. Hence~\eqref{eq:3.5} and~\eqref{eq:4.9} give $G_h^n\le \mathcal G_T+C_{\rm inv}h^{-1}\mathcal D_h(T)<M_*$, contradicting the choice of~$n$. Consequently, $G_h^n<M_*\le\kappa L_h$ whenever $t_n\le T$, so $T_h(\kappa)>T$ for all sufficiently small~$h$. Letting $T\uparrow T^*$ gives $\liminf_{h\downarrow0}T_h(\kappa)\ge T^*$.

For the upper bound, recall from~\eqref{eq:4.7} that $s_h<T^*$ and $G(s_h)=M_h$. Let $n_h$ be the last index with $t_{n_h}\le s_h$ and suppose, for contradiction, that $G_h^n<\kappa L_h$ for all $0\le n\le n_h$. Then condition~\eqref{eq:4.10} allows us to apply Lemma~\ref{Lemma:3.1} with $M_*=\kappa L_h$ and $M_2=\mathcal G_h(s_h)$. Applying the value estimate~\eqref{eq:3.5} at $t_{n_h}$ and extending $U_h^{n_h}$ by zero, we obtain the following bound from the definitions in~\eqref{eq:3.3} and~\eqref{eq:4.8}, since $0\le s_h-t_{n_h}\le\dt$:
\begin{equation}\label{eq:4.14}
 \norm{\overline U_h^{n_h}-u(s_h)}_{\infty,\Omega}
 \le \mathrm{Err}_h^{\rm ref}:=\mathcal D_h(s_h)+\varepsilon_h(s_h)+\beta_h(s_h)
                      +\omega_{\mathrm{time}}(\dt;s_h).
\end{equation}
To obtain a contradiction, choose $x_h\in\overline\Omega$ with
$|\nabla u(x_h,s_h)|=M_h$, and apply Lemma~\ref{Lemma:4.4} to
$q=\nabla u(x_h,s_h)$ to choose $e_h$. By \eqref{eq:4.5} and \eqref{eq:4.11}, the directional derivative
$\nabla u(x_h+se_h,s_h)\cdot e_h$ has constant sign and absolute value
at least $3a_\Omega M_h/4$ for $0\le s\le r_h$. Hence
\begin{equation*}
 |u(x_h+r_he_h,s_h)-u(x_h,s_h)|\ge\frac{3}{4}a_\Omega M_hr_h.
\end{equation*}
Combining~\eqref{eq:4.2} and~\eqref{eq:4.14} with this estimate gives
\begin{equation*}
 \begin{aligned}
 G_h^{n_h}&\ge
 \frac{|\overline U_h^{n_h}(x_h+r_he_h)-\overline U_h^{n_h}(x_h)|}{r_h}\\
 &\ge
 \frac{|u(x_h+r_he_h,s_h)-u(x_h,s_h)|
       -2\norm{\overline U_h^{n_h}-u(s_h)}_{\infty,\Omega}}{r_h}
 \ge\frac{3}{4}a_\Omega M_h-\frac{2\mathrm{Err}_h^{\rm ref}}{r_h}.
 \end{aligned}
\end{equation*}
By \eqref{eq:4.12}, the right-hand side is at least
$a_\Omega M_h/2>\kappa L_h$ for small $h$, a contradiction.
Hence $T_h(\kappa)\le t_{n_h}\le s_h<T^*$. Together with the lower bound,
this proves convergence.
Finally, by~\eqref{eq:4.1},
$T_h(\kappa_-)\le T_h(\kappa)\le T_h(\kappa_+)$ 
for all $\kappa\in[\kappa_-,\kappa_+]$. 
If the hypotheses hold at both endpoints, their convergence to~$T^*$ therefore gives uniform convergence~\eqref{eq:4.13}.
 \end{proof}

Regularity assumption~\hyperlink{regularity-R}{$\mathrm{(R)}$} gives explicit bounds for the two moduli of the exact solution and for the residual of the Ritz reference function~\eqref{eq:3.16}, all in terms of the exact gradient bound.
Applying Theorem~\ref{Theorem:4.5} then gives the following conditions on
$h$, $\dt$, and the prescribed threshold scale.
\begin{theorem}\label{Theorem:4.6}
Let $p>2$ and let $\Omega\subset\R^N$, $N\in\{1,2,3\}$, be a fixed bounded
$C^2$ domain.  Assume the domain, mesh, initialization, and exact-solve
hypotheses of Theorem~$\mathrm{\ref{Theorem:3.5}}$.
Assume that, for some~$T_{\mathrm{end}}>T^*$ independent of~$h$, each time grid specified in Section~$\mathrm{\ref{Section:2}}$ covers~$[0,T_{\mathrm{end}}]$.
Let $u$ be a classical solution of~\eqref{eq:1.1} satisfying regularity
assumption~\hyperlink{regularity-R}{$\mathrm{(R)}$} and \eqref{eq:1.2}.
Suppose that, as $h\downarrow0$,
\begin{equation}\label{eq:4.15}
 L_h\to\infty,\qquad hL_h^{3(p-1)}\to0,\qquad
              \dt h^{-1}L_h^{p-1}\to0.
\end{equation}
Then $T_h(\kappa)\to T^*$ for every fixed $\kappa>0$.
Moreover, for every $0<\kappa_-\le\kappa_+<\infty$ the convergence is uniform in the sense of~\eqref{eq:4.13}.
In particular, all these gradient-threshold times are finite for sufficiently
small $h$.
\end{theorem}
\begin{proof}
We verify the hypotheses of Theorem~\ref{Theorem:4.5} with the reference
 function \eqref{eq:3.16}. Estimate~\eqref{eq:3.17}
 gives uniformly bounded initial gradients, and
Lemma~\ref{Lemma:3.3} gives \eqref{eq:3.19}.
On each fixed interval $[0,T]$ with $T<T^*$, \eqref{eq:4.15} implies $\dt=o(h)$, so
$\mathcal G_h(T)\le \mathcal G_T$ and $\mathcal D_h(T)\le C_T h$ with constants independent
of $h$. This verifies condition~\textup{(i)}.
For condition~\textup{(ii)}, fix $\kappa>0$ and $C_*>2/a_\Omega$, and set $M_h=C_*\kappa L_h$ and
$s_h=s(M_h)$. By \eqref{eq:4.7}, $M(s_h)=M_h$,
so \eqref{eq:3.19} gives $\mathcal G_h(s_h)\le C_EM_h$. Thus, \eqref{eq:4.15} implies
\eqref{eq:2.30} with $M_1=\kappa L_h+C_EM_h$ and $M_2=C_EM_h$.
 Lemma~\ref{Lemma:2.8} therefore yields \eqref{eq:4.10}, verifying condition~\textup{(ii)}.
Under regularity assumption~\hyperlink{regularity-R}{$\mathrm{(R)}$}, the definitions in~\eqref{eq:4.8}, together with~\eqref{eq:3.11} and~\eqref{eq:3.12}, give
\begin{equation}\label{eq:4.16}
 \omega_{\mathrm{time}}(\delta;T)\le CM(T)^p\delta,\qquad
 \omega_{\mathrm{grad}}(r;T)\le C_{\rm reg}M(T)^p r.
\end{equation}
For condition~\textup{(iii)}, choose
\begin{equation}\label{eq:4.17}
 c_r=\min\left\{\ell_0,\frac{a_\Omega}{4(1+C_{\rm reg})}\right\}>0,
 \qquad r_h=c_rM_h^{-(p-1)}.
\end{equation}
By \eqref{eq:4.16} and \eqref{eq:4.17},
$4\omega_{\mathrm{grad}}(r_h;s_h)\le 4C_{\rm reg}M_h^pr_h\le a_\Omega M_h$,
 which verifies the first requirement in condition~\textup{(iii)}.
It remains to verify the second condition in~\textup{(iii)}. By~\eqref{eq:4.16},
\begin{equation*}
 \omega_{\mathrm{time}}(\dt;s_h)\le C\dt M_h^p\le C\dt M_h^{2p-2},
\end{equation*}
since $p>2$ and $M_h\ge1$ for sufficiently small $h$.
Thus, this error due to the time shift is absorbed into the right-hand side of the estimate for~$\mathcal D_h$ in~\eqref{eq:3.19}. Likewise,
$\varepsilon_h(s_h)$ and $\beta_h(s_h)$ are
absorbed into $ChM_h^{2p-1}$. Combining these bounds with~\eqref{eq:3.19}, we obtain
\begin{equation*}
 \mathrm{Err}_h^{\rm ref}
 =\mathcal D_h(s_h)+\varepsilon_h(s_h)+\beta_h(s_h)+\omega_{\mathrm{time}}(\dt;s_h)
 \le C(hM_h^{2p-1}+\dt M_h^{2p-2}).
\end{equation*}
Now, since $M_h=C_*\kappa L_h$, conditions~\eqref{eq:4.15} imply
\begin{equation*}
 hM_h^{3p-3}\to0,\qquad
 \dt M_h^{3p-4}
 =\bigl(\dt h^{-1}M_h^{p-1}\bigr)
  \bigl(hM_h^{3p-3}\bigr)M_h^{-p}\to0.
\end{equation*}
Combining these limits with the preceding error estimate and the choice of $r_h$ in~\eqref{eq:4.17}, we obtain
\begin{equation*}
 \frac{\mathrm{Err}_h^{\rm ref}}{M_hr_h}
 =c_r^{-1}\mathrm{Err}_h^{\rm ref}M_h^{p-2}
 \le C(hM_h^{3p-3}+\dt M_h^{3p-4})\to0.
\end{equation*}
This verifies~\eqref{eq:4.12}, so all the hypotheses of Theorem~\ref{Theorem:4.5} hold. Hence $T_h(\kappa)\to T^*$ for every fixed $\kappa>0$. Finally, applying that theorem at $\kappa_-$ and $\kappa_+$ gives the uniform conclusion~\eqref{eq:4.13}.
 \end{proof}

Substituting~$L_h=h^{-\sigma}$ into~\eqref{eq:4.15} gives the following corollary.
\begin{corollary}\label{Corollary:4.7}
Under the remaining hypotheses of Theorem~$\mathrm{\ref{Theorem:4.6}}$,
the threshold scale $L_h=h^{-\sigma}$ satisfies \eqref{eq:4.15}, and hence the conclusions of Theorem~$\mathrm{\ref{Theorem:4.6}}$ hold, if
\begin{equation*}
 0<\sigma<\frac{1}{3(p-1)},\qquad
 \dt=o(h^{1+\sigma(p-1)}).
\end{equation*}
\end{corollary}

\begin{remark}\label{Remark:4.8}
For inexact steps with exactly imposed numerical boundary values,
Theorem~\ref{Theorem:4.5} remains valid if every occurrence of
$\mathcal D_h(T)$ in its hypotheses is replaced by
\begin{equation*}
 \mathcal D_h(T)+\max_{t_n\le T}{\Sigma}_n,
 \qquad {\Sigma}_n=\sum_{k<n}\norm{\rho^{k+1}}_{\ell^\infty},
\end{equation*}
since the first inequality in~\eqref{eq:3.29} supplies this enlarged error bound. Thus it suffices to require
\begin{equation}\label{eq:4.18}
 \max_{t_n\le T}{\Sigma}_n=O_T(h)\quad\text{for every fixed }T<T^*,
 \qquad \frac{{\Sigma}_{n_h}}{M_hr_h}\to0,
\end{equation}
where $n_h$ is the last index with $t_{n_h}\le s_h$.
For the choice~\eqref{eq:4.17}, the residual tolerance~\eqref{eq:3.30} with fixed~$\eta>0$ ensures~\eqref{eq:4.18}; in particular, \eqref{eq:3.31} gives 
\begin{equation*}
 \frac{{\Sigma}_{n_h}}{M_hr_h}
 \le\frac{\eta T^*}{c_r}hM_h^{p-2}
 =\frac{\eta T^*}{c_r}(C_*\kappa)^{p-2}
       (hL_h^{3(p-1)})L_h^{-(2p-1)}\to0
\end{equation*}
under \eqref{eq:4.15}. Consequently, the conclusions of
Theorem~\ref{Theorem:4.6} and Corollary~\ref{Corollary:4.7} also hold
for these inexact solves. 
\end{remark}

\section{Numerical experiments}\label{Section:5}

We examine how mesh-dependent thresholds approximate the gradient blow-up time on inner approximations of the unit disk. To this end, 
we choose an initial value for which Appendix~\ref{Section:A} proves finite-time gradient blow-up.
We then compare the gradient-threshold times~$T_h(\kappa)$ for different~$\kappa$ as the mesh is refined. To compute these times, we solve the linear systems~\eqref{eq:2.14} by BiCGSTAB with a diagonal preconditioner~\cite{MR1149111}.

\paragraph{Problem and meshes.}
We take~$\Omega=\{x\in\R^2\colon |x|<1\}$, $p=5/2$, and
\begin{equation}\label{eq:5.1}
 u_0(x)=\begin{cases}
 13\exp\!\left(1-\dfrac1{1-4|x-(0.3,0)|^2}\right),
       &\text{for }|x-(0.3,0)|<0.5,\\
 0,&\text{otherwise}.
 \end{cases}
\end{equation}
For this choice, Corollary~\ref{Corollary:A.3} verifies the sufficient large-data condition of Proposition~\ref{Proposition:A.2}, thereby establishing finite-time gradient blow-up of the exact solution. Since the initial value is also smooth and compactly supported, Proposition~\ref{Proposition:3.2} gives the regularity assumption~\hyperlink{regularity-R}{(R)} required in Theorems~\ref{Theorem:3.5} and~\ref{Theorem:4.6}.

To approximate the disk from inside, we first triangulate the rhombus~$D\supset\overline\Omega$ of side length~$4$, centered at~$(0.3,0)$ with sides parallel to~$(1,0)$ and~$(1/2,\sqrt{3}/2)$, by equilateral triangles of side length~$h=1/J$, where~$J$ is a positive integer. From this background mesh we retain precisely the triangles whose three vertices lie in~$\Omega$ and take~$\Omega_h$ to be the interior of their union. By convexity, each retained triangle lies in~$\Omega$. Moreover, if~$x\in\Omega$ lies farther than~$h$ from~$\partial\Omega$, every background triangle containing~$x$ lies in~$\overline{B(x,h)}\subset\Omega$ and is retained. 
Hence~$\Omega_h$ satisfies~\eqref{eq:2.1}.
The background mesh~$\widehat{\mathcal T}_h$ and the computational mesh~$\mathcal T_h$ both consist of equilateral triangles of side length~$h$, so both satisfy~\eqref{eq:2.18}. Since~$\mathcal T_h\subset\widehat{\mathcal T}_h$ and~$\widehat h=h$, condition~\eqref{eq:2.36} also holds.
On the resulting computational mesh, six triangles meet at each interior node~$i$. To write the discrete operators at~$i$, we number its neighboring nodes~$j_k$ for~$k=0,\ldots,5$ counterclockwise, with~$d_k=(\cos(k\pi/3),\sin(k\pi/3))^{\mathsf{T}}$ and~$x_{j_k}=x_i+hd_k$. With this notation, the mesh coefficients are
\begin{equation*}
 m_i=\frac{\sqrt3}{2}h^2,\qquad
 a_{i,j_k}=\frac{2}{3h^2},\qquad \ell_{i,j_k}=\frac{d_k}{3h},
\end{equation*}
so the discrete Laplacian and recovered gradient take the form
\begin{equation*}
 (\Delta_h W)_i=\frac{2}{3h^2}\sum_{k=0}^5(W_{j_k}-W_i),\qquad
 g_i(W)=\frac{1}{3h}\sum_{k=0}^5 W_{j_k}d_k.
\end{equation*}
These coefficients satisfy the sign condition~\eqref{eq:2.15} and~\eqref{eq:2.22} with~$c_D=1/2$.

\paragraph{Threshold detection.}
We choose $L_h=h^{-1/8}$ and $\kappa=50,75,100$ in
\eqref{eq:4.1}, with the gradient norm taken over $\Omega_h$.
We take~$\tau=\min\{2^{-16},h^2\}$. The upper bound~$2^{-16}$ prevents excessively large time steps on the coarser meshes, where~$h^2$ is larger. Thus the time step is~$2^{-16}$ for~$1/h=32,64,128,256$ and~$h^2$ for~$1/h=512,1024$.
For~$p=5/2$, these choices satisfy the scale conditions of Corollary~\ref{Corollary:4.7} with~$\sigma=1/8$.
Unless otherwise stated, the linear solver stops at each time step once the residual $\rho^{n+1}$
defined in Section~\ref{Subsection:3.3} satisfies~\eqref{eq:3.30} with~$\eta=10^{-3}$, that is,
\begin{equation*}
 \norm{\rho^{n+1}}_{\ell^\infty}
 \le10^{-3}h\tau,
\end{equation*}
and this criterion was met at every computed step.
Each computation has a final-time limit of $0.005$ and is terminated earlier
at the first computed time with $h^{1/8}G_h^n\ge150$.

\begin{figure}[!ht]
\centering\includegraphics[width=\linewidth]{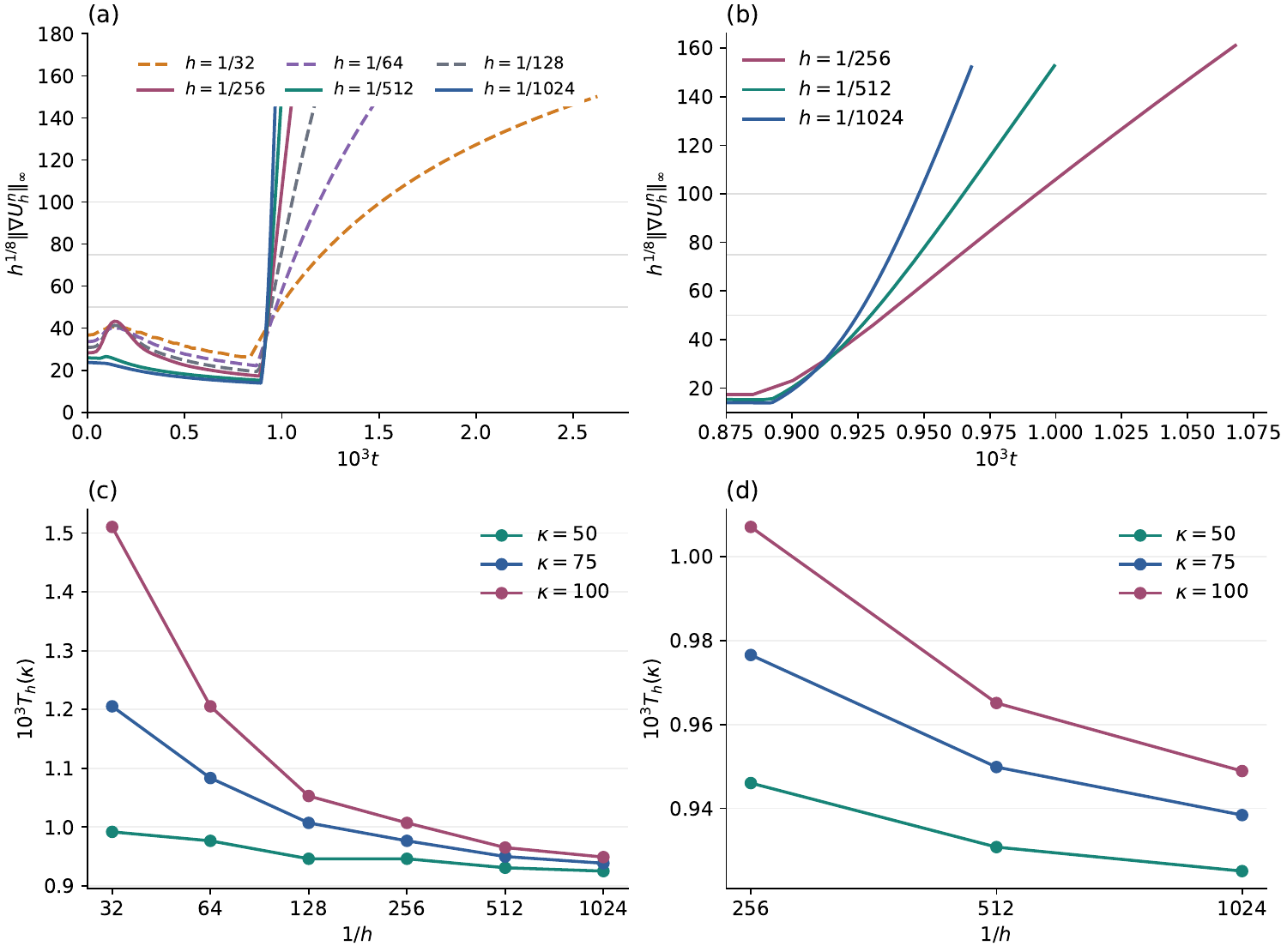}
\caption{Histories of $h^{1/8}G_h^n$ on all six meshes (a) and the three
finest meshes, $1/h=256,512,1024$ (b). In (a), dashed curves correspond to~$1/h=32,64,128$, and solid curves to~$1/h=256,512,1024$.
Panels (c) and (d) show the gradient-threshold times $T_h(\kappa)$ from
the computations in (a) and (b), respectively. The values in (c) are
reported in Table~\ref{tab:1}.
}\label{fig:1}\end{figure}

\paragraph{Results.}
Table~\ref{tab:1} reports the gradient-threshold times~$T_h(\kappa)$ defined in~\eqref{eq:4.1}.
Figure~\ref{fig:1}(a) shows that, after a short initial transient,~$h^{1/8}G_h^n$ decreases and then rises sharply on all six meshes. The enlarged view in panel (b)
shows that this rise occurs over a similar time interval on the three
finest meshes, $1/h=256,512,1024$.
On all six meshes, the solution values remain bounded throughout the computations: the maximum nodal value does not exceed its initial value, and the minimum nodal value stays above~$-6.80\times10^{-16}$.
Panels (c) and (d) compare the gradient-threshold times $T_h(\kappa)$ from the
same computations; the differences between the gradient-threshold times for $\kappa=50,75,100$
decrease under mesh refinement.
Refining the mesh from $h=1/512$ to $1/1024$ changes the gradient-threshold
times in Table~\ref{tab:1} by
0.61\%, 1.20\%, and 1.68\% for $\kappa=50,75,100$, respectively,
relative to the values at $h=1/512$.
To quantify the dependence on $\kappa$, we note that $T_h(100)-T_h(50)=2.38\times10^{-5}$ at $h=1/1024$.
This difference is~$2.58\%$ of~$T_h(50)$. 
Thus, even on the finest mesh, the estimated times vary by a few percent with the choice of threshold multiplier.

\begin{table}[ht]
\centering\small
\begin{tabular}{rrrrr}
\toprule
$1/h$ & $\tau$ & $10^3T_h(50)$ & $10^3T_h(75)$ & $10^3T_h(100)$\\
\midrule
 32 & $1.53\times10^{-5}$ & 0.991821 & 1.205444 & 1.510620\\
 64 & $1.53\times10^{-5}$ & 0.976562 & 1.083374 & 1.205444\\
 128 & $1.53\times10^{-5}$ & 0.946045 & 1.007080 & 1.052856\\
 256 & $1.53\times10^{-5}$ & 0.946045 & 0.976562 & 1.007080\\
 512 & $3.81\times10^{-6}$ & 0.930786 & 0.949860 & 0.965118\\
 1024 & $9.54\times10^{-7}$ & 0.925064 & 0.938416 & 0.948906\\
 \bottomrule
\end{tabular}
\caption{Gradient-threshold times $T_h(\kappa)$ for $\kappa=50,75,100$, with $\tau=\min\{2^{-16},h^2\}$.}
\label{tab:1}
\end{table}

\paragraph{Sensitivity to the time step.}
To assess the sensitivity of the gradient-threshold times to the time step, we perform one additional run at~$h=1/512$ with~$\tau=h^2/2$. This run has~$524$ time steps; the runs with~$\tau=h^2$ at~$h=1/512$ and~$h=1/1024$ reported in Table~\ref{tab:1} have~$262$ and~$1015$, respectively.
In this run, the maximum nodal value again stays at or below its initial value, and the minimum stays above~$-1.51\times10^{-64}$.
At $h=1/512$, halving the time step from $h^2$ to $h^2/2$ leaves
$T_h(50)$, $T_h(75)$, and $T_h(100)$ unchanged.
 To resolve the crossing times more finely than the time step, we also
estimate the time at which
 the numerical gradient reaches each threshold by linearly interpolating between
the last computed value of $G_h^n$ below $\kappa L_h$ and the first value
at or above it. Halving the time step changes these estimates by $2.95\times10^{-7}$,
$1.22\times10^{-7}$, and $2.98\times10^{-9}$ for $\kappa=50,75,100$,
respectively. Relative to the gradient-threshold times at~$h=1/512$ in Table~\ref{tab:1}, these changes are~$0.0317\%$,~$0.0128\%$, and~$0.000309\%$, respectively, showing little sensitivity to halving the time step in this test.

\paragraph{Solver tolerance and accumulated residuals.}
By~\eqref{eq:3.30}, a smaller~$\eta$ reduces the residual allowed at each time step and hence, by~\eqref{eq:3.31}, the bound on the accumulated residual.
To test how sensitive the computed threshold times are to this tolerance, we repeat the run at~$h=1/512$ and~$\tau=h^2$ with~$\eta$ reduced from~$10^{-3}$ to~$10^{-6}$. Tightening the tolerance leaves~$T_h(50)$,~$T_h(75)$, and~$T_h(100)$ unchanged and changes each linearly interpolated crossing-time estimate by less than~$10^{-15}$. Thus, in this test, tightening the residual tolerance changes the estimated crossing times far less than halving the time step does.

To assess the contribution of the linear solver residuals to the error bound~\eqref{eq:3.29}, Table~\ref{tab:2} reports the accumulated residual~$\Sigma_n$ at~$t_n=T_h(100)$. The column~$\Sigma_n/h$ compares this residual with the~$O(h)$ scale of the discretization error in that bound.
To account also for the prescribed tolerance, the table further reports~${\Sigma}_n/(\eta h)$; by~\eqref{eq:3.31}, this normalized ratio is at most~$t_n$. At~$h=1/512$ and~$\tau=h^2$, the ratios for~$\eta=10^{-3}$ and~$\eta=10^{-6}$ both lie below the common upper bound~$T_h(100)$.

\begin{table}[!htbp]
\centering\small
\begin{tabular}{rrrrrrr}
\toprule
$1/h$ & $\tau/h^2$ & $\eta$ & ${\Sigma}_n$ & ${\Sigma}_n/h$ & ${\Sigma}_n/(\eta h)$
  & $(100L_h)^{p-2}{\Sigma}_n$\\
\midrule
512 & $1$ & $10^{-3}$ & $8.0568\times10^{-11}$ & $4.1251\times10^{-8}$ & $4.1251\times10^{-5}$
                         & $1.1898\times10^{-9}$\\
512 & $1$ & $10^{-6}$ & $1.8252\times10^{-12}$ & $9.3453\times10^{-10}$ & $9.3453\times10^{-4}$
                         & $2.6956\times10^{-11}$\\
512 & $1/2$ & $10^{-3}$ & $7.0791\times10^{-11}$ & $3.6245\times10^{-8}$ & $3.6245\times10^{-5}$
                         & $1.0455\times10^{-9}$\\
1024 & $1$ & $10^{-3}$ & $3.6885\times10^{-11}$ & $3.7771\times10^{-8}$ & $3.7771\times10^{-5}$
                         & $5.6885\times10^{-10}$\\
\bottomrule
\end{tabular}
\caption{Accumulated solver residuals and their scaled values at~$t_n=T_h(100)$. The coefficient~$\eta$ specifies the stopping bound~\eqref{eq:3.30}.}
\label{tab:2}
\end{table}

The analysis of gradient-threshold times further requires the accumulated residual to be negligible relative to the spatial increment~$M_hr_h$ used to detect a large gradient; see Remark~\ref{Remark:4.8}. For the choice~$r_h=c_rM_h^{-(p-1)}$ in~\eqref{eq:4.17}, this comparison takes the form
\begin{equation*}
 \frac{{\Sigma}_n}{M_hr_h}=c_r^{-1}M_h^{p-2}{\Sigma}_n.
\end{equation*}
The last column therefore reports~$(100L_h)^{p-2}{\Sigma}_n$, using the detection threshold~$100L_h$ in place of the level~$M_h=100C_*L_h$ used in the proof. The table also uses the computed crossing time~$T_h(100)$, whereas the proof uses the last time level not exceeding~$s(M_h)$. At this crossing time, the scaled residual~$(100L_h)^{p-2}{\Sigma}_n$ is of order~$10^{-9}$ for~$\eta=10^{-3}$. For~$h=1/512$ and~$\tau=h^2$, it falls to~$2.6956\times10^{-11}$ when we tighten the tolerance to~$\eta=10^{-6}$.

\section*{Data availability}
The source code, mesh data, and numerical outputs needed to reproduce the computations, tables, and figures in Section~\ref{Section:5} are available from the corresponding author upon reasonable request.

\section*{Use of AI}
The authors used ChatGPT and Codex (OpenAI), with the GPT-6 Astra and GPT-6.1 Sol models, to assist in editing the text, developing and checking mathematical arguments, and writing code for the numerical experiments. They reviewed and verified all content, including the mathematical statements and proofs, and take full responsibility for the final manuscript.

\section*{Acknowledgements}
This work was partially supported by JSPS KAKENHI (grant number: 23K03179).

\appendix
\section{Regularity estimates and the disk example}\label{Section:A}
In this appendix, we prove the regularity estimates in Proposition~\ref{Proposition:3.2}. We also derive a sufficient condition for finite-time gradient blow-up and show that the initial data for the disk example in Section~\ref{Section:5} satisfy it.

\subsection{Maximum principle and proof of Proposition~\ref{Proposition:3.2}}
\label{Subsection:A.1}
We use the maximum principle to show that $u$ and $\partial_tu$ are uniformly bounded on~$\overline\Omega\times[0,T^*)$. In particular, we use the bound on~$\partial_tu$ to control the rescaled time derivative in the proof of Proposition~\ref{Proposition:3.2}.

\begin{proposition}\label{Proposition:A.1}
Let $u$ be a classical solution of \eqref{eq:1.1} on every compact interval
before $T^*$, with $\nabla u$ and $\partial_tu$ continuous up to $t=0$
and the boundary. Then, for $x\in\overline\Omega$ and $0\le t<T^*$,
\begin{equation}\label{eq:A.1}
 \min\{0,\min_{\overline\Omega}u_0\}\le u(x,t)
 \le\max\{0,\max_{\overline\Omega}u_0\},\qquad
 \norm{\partial_tu(t)}_{\infty,\Omega}
 \le\norm{\partial_tu(0)}_{\infty,\Omega}.
\end{equation}
\end{proposition}
\begin{proof}
The comparison principle \cite{PorreS2020}*{Proposition~3.3}
(see also \cite{Fried1964}*{Chapter~2, Section~6, Theorem~16})
gives the first bounds in \eqref{eq:A.1}. To estimate $\partial_tu$, we apply it to
$u(t+\varepsilon)$ and
$u(t)\pm\norm{u(\varepsilon)-u_0}_{\infty,\Omega}$ to obtain
\begin{equation*}
 \norm{u(t+\varepsilon)-u(t)}_{\infty,\Omega}
 \le\norm{u(\varepsilon)-u_0}_{\infty,\Omega}
\end{equation*}
for $0\le t<T^*-\varepsilon$. 
Divide this inequality by $\varepsilon$ and let
$\varepsilon\downarrow0$ to obtain the second bound in \eqref{eq:A.1}.
 \end{proof}

Regularity assumption~\hyperlink{regularity-R}{(R)} also implies the
hypotheses of Proposition~\ref{Proposition:A.1}, so that \eqref{eq:A.1}
holds for the solutions considered in Section~\ref{Section:3}.
Indeed, \eqref{eq:3.9} and \eqref{eq:3.13}--\eqref{eq:3.14} ensure that
$\nabla u$ and $\partial_tu$ extend continuously to $t=0$.

We now establish the estimates in Proposition~\ref{Proposition:3.2}.

\begin{proof}[Proof of Proposition~\ref{Proposition:3.2}]
We first establish regularity on a fixed initial interval. At later
times, we rescale space, time, and amplitude using the gradient bound
$M(T)$, apply local regularity estimates with uniform constants, and
return to the original variables.

\noindent\textit{Step 1: initial regularity.}
Set $H(\xi)=|\xi|^p$, fix $0<\alpha<\min\{1,p-2\}$, and choose
$q>(N+2)/(1-\alpha)$. Then $H\in C^{2,\alpha}_{\rm loc}(\R^N)$.
On a fixed initial interval $[0,2t_*]\subset[0,T^*)$, boundedness of
$\nabla u$, the parabolic $W^{2,1}_q(\Omega\times(0,2t_*))$ estimate, and embedding give
$\nabla u\in C^{\alpha,\alpha/2}(\overline\Omega\times[0,2t_*])$; see
\cite{LadySU1968}*{Chapter~IV, Theorem~9.1; Chapter~II, Lemma~3.3}.
Define
\begin{equation*}
 v_0=\Delta u_0+H(\nabla u_0),\qquad
 v_1=\Delta v_0+DH(\nabla u_0)\cdot\nabla v_0.
\end{equation*}
Since $u_0\in C_0^5(\Omega)$, $v_0\in C^{2,\alpha}(\overline\Omega)$,
and $u_0,v_0,v_1$ vanish near $\partial\Omega$. The compatibility
conditions above and the Schauder theorem
\cite{LadySU1968}*{Chapter~IV, Theorem~5.2} give
$u\in C^{2+\alpha,1+\alpha/2}(\overline\Omega\times[0,2t_*])$ and $\partial_tu(0)=v_0$.
Time difference quotients then show that $v=\partial_tu$ satisfies
\begin{equation*}
 \partial_t v-\Delta v=DH(\nabla u)\cdot\nabla v,\qquad
 v|_{\partial\Omega}=0,\qquad v(0)=v_0.
\end{equation*}
The same Schauder theorem applied to this equation gives
$v\in C^{2+\alpha,1+\alpha/2}(\overline\Omega\times[0,2t_*])$, since $v_1=0$ on the boundary.
This proves~\eqref{eq:3.8} on the initial interval and gives the required initial regularity.
Proposition~\ref{Proposition:A.1}
now yields $\norm{v(t)}_{\infty,\Omega}\le A_0:=\norm{v_0}_{\infty,\Omega}$
for every $t<T^*$.

\noindent\textit{Step 2: rescaling.}
Fix a finite covering of $\partial\Omega$ by open coordinate neighborhoods
in $\R^N$ in which the boundary is a smooth graph.
By compactness, choose $r_0>0$ such that, for every $x_c\in\partial\Omega$,
the ball $B(x_c,r_0)$ is contained in one of these neighborhoods.
Choose a fixed $0<c_\mu\le1$ with $64c_\mu^2<t_*/2$ and $8c_\mu<r_0$.
For $T\ge t_*$, put $M=M(T)$ and $\mu=c_\mu M^{-(p-1)}$.
Since $M\ge1$, we have $8\mu\le8c_\mu<r_0$; hence every $B(x_c,8\mu)$,
$x_c\in\partial\Omega$, is contained in one of the fixed coordinate neighborhoods.

Fix $t_c\in[t_*,T]$.
We estimate derivatives at each point $x\in\overline\Omega$ by working
on a local spatial neighborhood and a backward time interval ending at $t_c$.
If $\operatorname{dist}(x,\partial\Omega)\ge2\mu$, take $x_c=x$ and use
the interior ball $B(x_c,2\mu)\subset\Omega$.
Otherwise, choose a nearest boundary point $x_c\in\partial\Omega$ and
use the boundary neighborhood $B(x_c,8\mu)\cap\Omega$.
In the first case, $x=x_c$ is the center of the smaller ball $B(x_c,\mu)$.
In the second, $|x-x_c|<2\mu$, so $x$ lies in the smaller neighborhood
$B(x_c,4\mu)\cap\overline\Omega$.
Set $k=u(x_c,t_c)$ for an interior neighborhood and $k=0$ for a boundary
neighborhood. The rescaled functions are defined by
\begin{equation}\label{eq:A.2}
 \begin{aligned}
 \widetilde u(y,s)&:=\mu^{-(p-2)/(p-1)}
       \bigl(u(x_c+\mu y,t_c+\mu^2s)-k\bigr),\\
 \widetilde v(y,s)&:=\partial_s \widetilde u(y,s)
       =\mu^{p/(p-1)}v(x_c+\mu y,t_c+\mu^2s).
 \end{aligned}
\end{equation}
In the rescaled neighborhoods, we choose nested backward cylinders
$\mathcal Q_0\supset \mathcal Q_1\supset \mathcal Q_2\supset \mathcal Q_3$,
whose spatial balls are centered at $y=0$ and whose time intervals end at
$s=0$, with this endpoint included. We take outer and inner spatial
radii $2$ and $1$ in the interior case, and $8$ and $4$ in the boundary
case, choosing the intermediate radii with fixed positive gaps.
For each cylinder, the length of its time interval is the square of its spatial radius.
Thus the outer time interval has length at most $64$, and our choice of $c_\mu$ ensures
that all original time intervals lie in $[t_*/2,T]$.

By the choice of the inner radius, the rescaled point~$\bigl((x-x_c)/\mu,0\bigr)$ lies in~$\overline{\mathcal Q_3}$. To estimate derivatives there, we start on the largest cylinder~$\mathcal Q_0$. Using~\eqref{eq:1.1}, we obtain the rescaled equation and basic bounds
\begin{equation}\label{eq:A.3}
 \begin{cases}
 \partial_s \widetilde u-\Delta_y\widetilde u=H(D\widetilde u)
       \quad\text{in }\mathcal Q_0,\\[2pt]
 \norm{D\widetilde u}_{\infty,\mathcal Q_0}\le\mu^{1/(p-1)}M=c_\mu^{1/(p-1)},\\[2pt]
 \norm{\widetilde v}_{\infty,\mathcal Q_0}\le A_0\mu^{p/(p-1)}\le C,\\[2pt]
 \norm{\widetilde u}_{\infty,\mathcal Q_0}
 \le C\bigl(\mu^{1/(p-1)}M+A_0\mu^{p/(p-1)}\bigr)\le C.
 \end{cases}
\end{equation}
Under the rescaling~\eqref{eq:A.2}, the estimates~$|\nabla u|\le M$ and~$|\partial_tu|\le A_0$ give the bounds for~$D\widetilde u$ and~$\widetilde v$. In an interior neighborhood, these same estimates control the difference~$u(x_c+\mu y,t_c+\mu^2s)-u(x_c,t_c)$ by the mean value theorem, which gives the bound for~$\widetilde u$. In a boundary neighborhood, where~$k=0$, the corresponding bound follows from the zero boundary condition through~$|u(x,t)|\le M\operatorname{dist}(x,\partial\Omega)$. Since~$\mu^{1/(p-1)}M=c_\mu^{1/(p-1)}$ and~$\mu\le c_\mu\le1$, all three bounds in~\eqref{eq:A.3} are uniform in~$t_c$ and~$T$.

\noindent\textit{Step 3: estimates in the rescaled variables.}
To apply boundary regularity estimates to~\eqref{eq:A.3} with constants independent of~$\mu$, we first bound the rescaled boundary graph. In a fixed coordinate neighborhood containing~$x_c=(x_c',x_{c,N})$, write the boundary as~$x_N=\psi(x')$.
Since $x_{c,N}=\psi(x_c')$, the change of variables $x=x_c+\mu y$
gives the rescaled graph
$y_N=\psi_\mu(y'):=\mu^{-1}\bigl[\psi(x_c'+\mu y')-\psi(x_c')\bigr]$
and derivative bounds
\begin{equation*}
 D^j\psi_\mu=\mu^{j-1}(D^j\psi)(x_c'+\mu\,\cdot)
 \quad\mbox{for }j=1,2,3, \qquad
 [D^3\psi_\mu]_{\alpha,B'_8}
 \le\mu^{2+\alpha}[D^3\psi]_{\alpha,x_c'+B'_{8\mu}}.
\end{equation*}
Here $B'_r\subset\R^{N-1}$ is the ball of radius $r$ centered at the origin, and
$[\cdot]_{\alpha,E}$ is the spatial H\"older seminorm on a set~$E$.
These estimates give uniform $C^{3,\alpha}$
bounds for the flattening maps and their inverses. The transformed
operators are therefore uniformly elliptic with uniformly bounded
$C^{1,\alpha}$ coefficients.

We may therefore apply the local $W^{2,1}_q$ estimate
\cite{LadySU1968}*{Chapter~IV, Section~10} and embedding
\cite{LadySU1968}*{Chapter~II, Lemma~3.3} to \eqref{eq:A.3}, obtaining a uniform
 $C^{\alpha,\alpha/2}(\mathcal Q_1)$ bound for $D\widetilde u$. Thus $H(D\widetilde u)$ and $DH(D\widetilde u)$
have uniform H\"older bounds. Since $\widetilde u$ vanishes
on the part of the lateral boundary inherited from $\partial\Omega$,
the interior and boundary Schauder
estimates \cite{LadySU1968}*{Chapter~IV, Theorem~10.1} imply
\begin{equation}\label{eq:A.4}
 \norm{\widetilde u}_{C^{2+\alpha,1+\alpha/2}(\mathcal Q_2)}\le C.
\end{equation}
Next, for sufficiently small $\varepsilon>0$, define the backward time
difference quotients on $\mathcal Q_2$ by
\begin{equation*}
 w_\varepsilon(y,s)
 :=\frac{\widetilde u(y,s)-\widetilde u(y,s-\varepsilon)}{\varepsilon}
 =\int_0^1\widetilde v(y,s-\theta\varepsilon)\,\mathrm{d}\theta.
\end{equation*}
These quotients solve linear equations whose drifts are averages of $DH$
between the corresponding gradients, so the H\"older bounds on
$\mathcal Q_1$ established above give uniform
$C^{\alpha,\alpha/2}$ bounds for these drifts.
The integral representation and \eqref{eq:A.3} give
$|w_\varepsilon|\le A_0\mu^{p/(p-1)}$, uniformly in $\varepsilon$.
Moreover, $w_\varepsilon=0$ on the boundary corresponding to
$\partial\Omega$, since $\widetilde u=0$ there.
Applying the same Schauder estimate and letting $\varepsilon\downarrow0$
therefore gives
\begin{equation}\label{eq:A.5}
 \partial_s \widetilde v-\Delta \widetilde v=DH(D\widetilde u)\cdot D\widetilde v,\qquad
 \norm{\widetilde v}_{C^{2+\alpha,1+\alpha/2}(\mathcal Q_3)}\le C\mu^{p/(p-1)}.
\end{equation}
The zero boundary condition passes to the limit, so $\widetilde v=0$
on this boundary as well.

\noindent\textit{Step 4: scaling back.}
Under the change of variables~\eqref{eq:A.2}, the equation for~$\widetilde v$ in~\eqref{eq:A.5} gives~\eqref{eq:3.8} for~$t_*\le t\le T$.
Applying the same change of variables to the estimates in~\eqref{eq:A.4}--\eqref{eq:A.5},
we obtain bounds of orders $\mu^{-p/(p-1)}$, $\mu^{-1}$, and $\mu^{-2}$ for
$D^2u$, $\nabla v$, and $(D^2v,\partial_tv)$, respectively. Since $\mu=c_\mu M^{-(p-1)}$,
taking the supremum over the covering neighborhoods gives
$L^\infty(\Omega)$ bounds of orders $M^p$, $M^{p-1}$, and $M^{2p-2}$,
respectively, with constants independent of $t_c$ and $T$.
Since~$M\ge1$, the bound of order~$M^{2p-2}$ for~$D^2v$ also gives the order~$M^{2p-1}$ required in~\eqref{eq:3.10}. Combining these estimates with the fixed initial-interval bounds proves~\eqref{eq:3.9}--\eqref{eq:3.10} for all~$T<T^*$.
 \end{proof}

\subsection{A large-data criterion and the disk example}

Following the argument in~\cite{Soupl2002a}, we give
an explicit condition on the initial data ensuring finite-time gradient
blow-up.

\begin{proposition}
\label{Proposition:A.2}
Let $\Omega\subset\R^N$, $N\ge1$, be a smooth bounded connected domain, let $p>2$,
and let $u$ be the classical solution of \eqref{eq:1.1} with
$u_0\in C^1(\overline\Omega)$, $u_0\ge0$, and $u_0=0$ on
$\partial\Omega$. Let $\lambda_1>0$ and $\phi_1$ be the first Dirichlet
eigenvalue and eigenfunction of $-\Delta$, normalized by
$\phi_1>0$ in $\Omega$ and $\norm{\phi_1}_{L^1(\Omega)}=1$. Define
\begin{equation}\label{eq:A.6}
 B_{\Omega,p}:=\frac{1}{\lambda_1}
 \int_\Omega |\nabla\phi_1|^{\frac{p}{p-1}}
                   \phi_1^{-\frac{1}{p-1}}\, \mathrm{d}x,
 \qquad z_0:=\int_\Omega u_0\phi_1\, \mathrm{d}x.
\end{equation}
Then $0<B_{\Omega,p}<\infty$. If $z_0>B_{\Omega,p}$, the solution
exhibits finite-time gradient blow-up, and
\begin{equation}\label{eq:A.7}
 T^*\le\frac{1}{(p-1)\lambda_1}
 \log\!\left(\frac{1}{1-(B_{\Omega,p}/z_0)^{p-1}}\right).
\end{equation}
\end{proposition}
\begin{proof}
By~\cite{Soupl2002a}, $B_{\Omega,p}\in(0,\infty)$, and finite-time
breakdown is gradient blow-up.
It remains to prove~\eqref{eq:A.7}.
Assume $z_0>B_{\Omega,p}$ and set
\[ z(t)=\int_\Omega u(t)\phi_1\,\mathrm{d}x>0. \]
Integration by parts and the H\"older inequality give
\(
 \partial_t z\ge\lambda_1B_{\Omega,p}^{1-p}z^p-\lambda_1z.
\)
Multiplication by $(1-p)z^{-p}$ and integration give
\begin{equation*}
 z(t)^{1-p}\le B_{\Omega,p}^{1-p}
  +\bigl(z_0^{1-p}-B_{\Omega,p}^{1-p}\bigr)
                   e^{(p-1)\lambda_1t}.
\end{equation*}
The right-hand side vanishes at the time in~\eqref{eq:A.7},
whereas $z(t)^{1-p}>0$ for $t<T^*$.
This proves~\eqref{eq:A.7}.
 \end{proof}

For the initial data in~\eqref{eq:5.1}, it remains to compare
the weighted initial value $z_0$ with $B_{\Omega,p}$. We do so
by bounding the two integrals in~\eqref{eq:A.6} from opposite sides.

\begin{corollary}
\label{Corollary:A.3}
Let $\Omega=\{x\in\R^2\colon |x|<1\}$, $p=5/2$, and let $u_0$ be the
initial data in~\eqref{eq:5.1}.
The quantities in \eqref{eq:A.6} satisfy
\begin{equation}\label{eq:A.8}
 B_{\Omega,5/2}<2.312560<2.412749<z_0.
\end{equation}
 Consequently, the solution
 exhibits finite-time
gradient blow-up.
\end{corollary}
\begin{proof}
To estimate the two integrals in~\eqref{eq:A.6}, we use the representation of the first Dirichlet eigenfunction on the unit disk in terms of the Bessel function~$J_0$.
For the properties of Bessel functions used below, see
\cite{AbramS1964}*{Chapter~9}.
We reduce both integrals to one dimension and use the change of variables
$r=1-s^3$ to rewrite the integral for $B_{\Omega,5/2}$ with a continuous
integrand. We then evaluate rigorous upper and lower sums by interval
arithmetic.

Write $J_0,J_1$ for the Bessel functions of the first kind of orders
zero and one, respectively, and let
$j_{0,1}$ be the first positive zero of $J_0$.
By the Bessel equation,
$J_0(j_{0,1}|x|)$ satisfies
$-\Delta J_0(j_{0,1}|x|)=j_{0,1}^2J_0(j_{0,1}|x|)$ on the unit disk.
Since it vanishes on the boundary and is positive in the interior,
the first Dirichlet eigenpair, with the normalization
of Proposition~\ref{Proposition:A.2}, is
\begin{equation*}
 \lambda_1=j_{0,1}^2,\qquad \phi_1(x)=c_{\phi}J_0(j_{0,1}|x|),\qquad
 c_{\phi}=\frac{j_{0,1}}{2\pi J_1(j_{0,1})},
\end{equation*}
since the Bessel recurrence relations give
\( (rJ_1(j_{0,1}r))' 
 =j_{0,1}rJ_0(j_{0,1}r).
\)
Substituting this eigenpair in \eqref{eq:A.6} and
integrating in polar coordinates give
\begin{equation*}
 B_{\Omega,5/2}
 =\frac{j_{0,1}^{2/3}}{J_1(j_{0,1})}
   \int_0^1 rJ_1(j_{0,1}r)^{\frac{5}{3}}J_0(j_{0,1}r)^{-\frac{2}{3}}\,\mathrm{d}r.
\end{equation*}
To remove the singularity at~$r=1$, put~$r=1-s^3$. For~$0<s<1$, define
\begin{equation}\label{eq:A.9}
 F(s)=\frac{3j_{0,1}^{2/3}}{J_1(j_{0,1})}(1-s^3)
 \frac{J_1(j_{0,1}(1-s^3))^{5/3}}{[J_0(j_{0,1}(1-s^3))/s^3]^{2/3}},\qquad
 B_{\Omega,5/2}=\int_0^1F(s)\,\mathrm{d}s.
\end{equation}
Since $j_{0,1}$ is the first positive zero of $J_0$ and is simple,
$J_1(j_{0,1})=-J_0'(j_{0,1})>0$. The definition of the derivative
and $J_0(j_{0,1})=0$ give
\begin{equation*}
 \lim_{s\downarrow0}\frac{J_0(j_{0,1}(1-s^3))}{s^3}
 =j_{0,1}J_1(j_{0,1})>0.
\end{equation*}
Thus the denominator $[J_0(j_{0,1}(1-s^3))/s^3]^{2/3}$ in
\eqref{eq:A.9} has a positive limit, and $F$ extends continuously
with $F(0)=3$ and $F(1)=0$.

We next compute~$z_0$, using the fact that~$u_0$ is supported in the closed disk of radius~$1/2$ centered at~$(0.3,0)$, which lies inside~$\Omega$.
Let $r=|x-(0.3,0)|$ denote the distance from the center of the support.
Integrating Graf's addition theorem in the angular variable gives
\begin{equation*}
 \int_0^{2\pi}
 J_0\!\left(j_{0,1}\sqrt{a_0^2+r^2+2a_0r\cos\theta}\right)\,\mathrm{d}\theta
 =2\pi J_0(j_{0,1}a_0)J_0(j_{0,1}r).
\end{equation*}
Using this identity with $a_0=0.3$ in \eqref{eq:A.6}, we obtain
\begin{equation}\label{eq:A.10}
 z_0=\frac{13j_{0,1}J_0(0.3j_{0,1})}{J_1(j_{0,1})}
 \int_0^{1/2}r\exp\!\left(-\frac{4r^2}{1-4r^2}\right)
                         J_0(j_{0,1}r)\,\mathrm{d}r.
\end{equation}

We now bound the two one-dimensional integrals using fixed-point interval arithmetic. Every
interval endpoint is an integer multiple of $10^{-60}$, and products,
quotients, and cube roots are rounded outward. For $\iota=0,1$, use
the series expansion
\begin{equation*}
 J_\iota(x)=\left(\frac{x}{2}\right)^\iota
 \sum_{k=0}^{\infty}\frac{(-1)^k(x^2/4)^k}{k!(k+\iota)!}.
\end{equation*}
For $0\le x\le2.405$, the partial sums of these series
ending at $k=41$ and $k=40$ are lower and upper bounds, respectively,
because the omitted tail alternates with decreasing terms.
For $j_{0,1}^{-}=2.404825557695$ and $j_{0,1}^{+}=2.404825557696$, interval evaluation gives $J_0(j_{0,1}^{-})>0>J_0(j_{0,1}^{+})$.
Also, the series for $J_1$ gives
$J_1(x)\ge(x/2)(1-x^2/8)>0$ for $0<x<2.405$.
Hence $J_0$ is strictly decreasing there, and the endpoints above
satisfy $j_{0,1}^{-}<j_{0,1}<j_{0,1}^{+}$.
For the upper bound on~$B_{\Omega,5/2}$, we partition~$[0,1]$ into~$4096$ equal subintervals in~$s$. We then use interval arithmetic to obtain an upper bound for~$F$ in~\eqref{eq:A.9} on each subinterval and sum these bounds multiplied by the subinterval lengths.
For each subinterval, let~$r_-$ and~$r_+$ be lower and upper bounds for~$r=1-s^3$. Using the series bounds above, we obtain an upper bound for~$J_1(j_{0,1}r)$ by interval arithmetic on~$[j_{0,1}^{-}r_-,j_{0,1}^{+}r_+]$. Raising this bound to the power~$5/3$ and multiplying by~$r_+$ then bounds the numerator from above. We also bound the prefactor~$3j_{0,1}^{2/3}/J_1(j_{0,1})$ from above using the bounds on~$j_{0,1}$ and the series bounds.
To bound the denominator from below on the first interval, note that $s^3=1-r$, and use $J_0'=-J_1$ and $J_0(j_{0,1})=0$ to write
\begin{equation*}
 \frac{J_0(j_{0,1}r)}{1-r}
 =\frac{1}{1-r}\int_{j_{0,1}r}^{j_{0,1}}J_1(y)\,\mathrm{d}y.
\end{equation*}
This bounds
$J_0(j_{0,1}r)/(1-r)$ below by $j_{0,1}$ times a lower bound for $J_1$ on
$[j_{0,1}r,j_{0,1}]$. On any later interval, let $\bar s$ be its right
endpoint. Monotonicity gives the
positive lower bound
\begin{equation*}
 \frac{J_0(j_{0,1}r)}{s^3}
 \ge\frac{J_0(j_{0,1}^{+}r_+)}{\bar s^3}.
\end{equation*}
We use this lower bound to control the denominator in~\eqref{eq:A.9}.
To obtain the lower bound on~$z_0$ in~\eqref{eq:A.10}, we partition~$[0,0.48]$ into~$4096$ equal subintervals in~$r$ and discard the nonnegative contribution on~$[0.48,1/2]$.
For the exponential factor in \eqref{eq:A.10}, set $\rho=4r^2/(1-4r^2)$.
For $0\le r\le0.48$, we have $0\le\rho/16<1$, so the alternating
Taylor series gives
\begin{equation*}
 0<\sum_{k=0}^{41}\frac{(-\rho/16)^k}{k!}
 \le \exp(-\rho/16).
\end{equation*}
Taking sixteenth powers gives a lower bound for
$\exp(-\rho)=\bigl(\exp(-\rho/16)\bigr)^{16}$.
Using this exponential bound and the series bound for~$J_0$, we evaluate both decreasing functions at the right endpoint of each subinterval to obtain the required lower bounds.
These calculations yield \eqref{eq:A.8}. Applying
Proposition~\ref{Proposition:A.2} now gives finite-time gradient blow-up.
 \end{proof}


\begin{bibdiv}
\begin{biblist}
\bib{AbramS1964}{book}{
    author={Abramowitz, Milton},
    author={Stegun, Irene A.},
    title={Handbook of mathematical functions with formulas, graphs, and mathematical tables},
    series={National Bureau of Standards Applied Mathematics Series},
    volume={No. 55},
    publisher={U. S. Government Printing Office, Washington, DC},
    date={1964},
    pages={xiv+1046},
    note={For sale by the Superintendent of Documents},
}

{\color{black}
\bib{AcostDR2002}{article}{
   author={Acosta, G.},
   author={Dur\'an, R. G.},
   author={Rossi, J. D.},
   title={An adaptive time step procedure for a parabolic problem with
   blow-up},
   journal={Computing},
   volume={68},
   date={2002},
   number={4},
   pages={343--373},
   issn={0010-485X},
   doi={10.1007/s00607-002-1449-x},
}
}

{\color{black}
\bib{AttouS2020}{article}{
   author={Attouchi, Amal},
   author={Souplet, Philippe},
   title={Gradient blow-up rates and sharp gradient estimates for diffusive
   Hamilton-Jacobi equations},
   journal={Calc. Var. Partial Differential Equations},
   volume={59},
   date={2020},
   number={5},
   pages={Paper No. 153, 28},
   issn={0944-2669},
   doi={10.1007/s00526-020-01810-9},
}
}

{\color{black}
\bib{BabaT1981}{article}{
    author={Baba, Kinji},
    author={Tabata, Masahisa},
    title={On a conservative upwind finite element scheme for convective diffusion equations},
    journal={RAIRO Anal. Num\'er.},
    volume={15},
    date={1981},
    number={1},
    pages={3--25},
    issn={0399-0516,0516-2777},
}
}

{\color{black}
\bib{BarleD2004}{article}{
    author={Barles, Guy},
    author={Da Lio, Francesca},
    title={On the generalized {D}irichlet problem for viscous {H}amilton-{J}acobi equations},
    journal={J. Math. Pures Appl. (9)},
    volume={83},
    date={2004},
    number={1},
    pages={53--75},
    doi={10.1016/S0021-7824(03)00070-9},
    url={https://doi.org/10.1016/S0021-7824(03)00070-9},
    issn={0021-7824},
}
}

{\color{black}
\bib{BergeK1988}{article}{
    author={Berger, Marsha},
    author={Kohn, Robert V.},
    title={A rescaling algorithm for the numerical calculation of blowing-up solutions},
    journal={Comm. Pure Appl. Math.},
    volume={41},
    date={1988},
    number={6},
    pages={841--863},
    doi={10.1002/cpa.3160410606},
    url={https://doi.org/10.1002/cpa.3160410606},
    issn={0010-3640},
}
}

{\color{black}
\bib{BrambDT1972}{article}{
   author={Bramble, James H.},
   author={Dupont, Todd},
   author={Thom\'ee, Vidar},
   title={Projection methods for Dirichlet's problem in approximating
   polygonal domains with boundary-value corrections},
   journal={Math. Comp.},
   volume={26},
   date={1972},
   pages={869--879},
   issn={0025-5718},
   doi={10.2307/2005869},
}
}

{\color{black}
\bib{BrennS2008}{book}{
   author={Brenner, Susanne C.},
   author={Scott, L. Ridgway},
   title={The mathematical theory of finite element methods},
   series={Texts in Applied Mathematics},
   volume={15},
   edition={3},
   publisher={Springer, New York},
   date={2008},
   pages={xviii+397},
   isbn={978-0-387-75933-3},
   doi={10.1007/978-0-387-75934-0},
}
}

{\color{black}
\bib{BrunnWZ2010}{article}{
   author={Brunner, Hermann},
   author={Wu, Xiaonan},
   author={Zhang, Jiwei},
   title={Computational solution of blow-up problems for semilinear
   parabolic PDEs on unbounded domains},
   journal={SIAM J. Sci. Comput.},
   volume={31},
   date={2009/10},
   number={6},
   pages={4478--4496},
   issn={1064-8275},
   doi={10.1137/090761367},
}
}

{\color{black}
\bib{CangiGKM2016}{article}{
   author={Cangiani, Andrea},
   author={Georgoulis, Emmanuil H.},
   author={Kyza, Irene},
   author={Metcalfe, Stephen},
   title={Adaptivity and blow-up detection for nonlinear evolution problems},
   journal={SIAM J. Sci. Comput.},
   volume={38},
   date={2016},
   number={6},
   pages={A3833--A3856},
   issn={1064-8275},
   doi={10.1137/16M106073X},
}
}

{\color{black}
\bib{Cho2013}{article}{
   author={Cho, Chien-Hong},
   title={On the computation of the numerical blow-up time},
   journal={Jpn. J. Ind. Appl. Math.},
   volume={30},
   date={2013},
   number={2},
   pages={331--349},
   issn={0916-7005},
   doi={10.1007/s13160-013-0101-9},
}
}

{\color{black}
\bib{ChoS2026}{article}{
   author={Cho, Chien-Hong},
   author={Sun, Hao-Wei},
   title={On the numerical approximation of a rescaling algorithm to
   nonlinear blow-up problems},
   journal={Numer. Math.},
   volume={158},
   date={2026},
   number={1},
   pages={103--141},
   issn={0029-599X},
   doi={10.1007/s00211-025-01510-6},
}
}

{\color{black}
\bib{DemloLSW2012}{article}{
   author={Demlow, A.},
   author={Leykekhman, D.},
   author={Schatz, A. H.},
   author={Wahlbin, L. B.},
   title={Best approximation property in the $W^{1}_{\infty}$ norm for
   finite element methods on graded meshes},
   journal={Math. Comp.},
   volume={81},
   date={2012},
   number={278},
   pages={743--764},
   issn={0025-5718},
   doi={10.1090/S0025-5718-2011-02546-9},
}
}

{\color{black}
\bib{DieniRS2024}{article}{
   author={Diening, Lars},
   author={Rolfes, Julian},
   author={Salgado, Abner J.},
   title={Pointwise gradient estimate of the Ritz projection},
   journal={SIAM J. Numer. Anal.},
   volume={62},
   date={2024},
   number={3},
   pages={1212--1225},
   issn={0036-1429},
   doi={10.1137/23M1571800},
}
}

{\color{black}
\bib{Evans1998}{book}{
    author={Evans, Lawrence C.},
    title={Partial differential equations},
    series={Graduate Studies in Mathematics},
    volume={19},
    publisher={American Mathematical Society, Providence, RI},
    date={1998},
    pages={xviii+662},
    isbn={0-8218-0772-2},
}
}

{\color{black}
\bib{Fried1964}{book}{
    author={Friedman, Avner},
    title={Partial differential equations of parabolic type},
    publisher={Prentice-Hall, Inc., Englewood Cliffs, N.J.},
    date={1964},
    pages={xiv+347},
}
}

{\color{black}
\bib{Fujit1966}{article}{
    author={Fujita, Hiroshi},
    title={On the blowing up of solutions of the {C}auchy problem for {$u_{t}=\Delta u+u^{1+\alpha }$}},
    journal={J. Fac. Sci. Univ. Tokyo Sect. I},
    volume={13},
    date={1966},
    pages={109--124 (1966)},
    issn={0040-8980},
}
}

{\color{black}
\bib{GaoQ2021}{article}{
   author={Gao, Huadong},
   author={Qiu, Weifeng},
   title={The pointwise stabilities of piecewise linear finite element
   method on non-obtuse tetrahedral meshes of nonconvex polyhedra},
   journal={J. Sci. Comput.},
   volume={87},
   date={2021},
   number={2},
   pages={Paper No. 53, 5},
   issn={0885-7474},
   doi={10.1007/s10915-021-01465-4},
}
}

{\color{black}
\bib{GilbaT1983}{book}{
    author={Gilbarg, David},
    author={Trudinger, Neil S.},
    title={Elliptic partial differential equations of second order},
    series={Grundlehren der Mathematischen Wissenschaften [Fundamental Principles of Mathematical Sciences]},
    volume={224},
    publisher={Springer-Verlag, Berlin},
    date={1983},
    pages={xiii+513},
    edition={2},
    doi={10.1007/978-3-642-61798-0},
    url={http://dx.doi.org/10.1007/978-3-642-61798-0},
    isbn={3-540-13025-X},
}
}

{\color{black}
\bib{GongLX2023}{article}{
   author={Gong, Wei},
   author={Liang, Dongdong},
   author={Xie, Xiaoping},
   title={Pointwise error estimates for linear finite element approximation
   to elliptic Dirichlet problems in smooth domains},
   journal={Adv. Comput. Math.},
   volume={49},
   date={2023},
   number={2},
   pages={Paper No. 15, 33},
   issn={1019-7168},
   doi={10.1007/s10444-023-10017-3},
}
}

{\color{black}
\bib{Grois2006}{article}{
   author={Groisman, Pablo},
   title={Totally discrete explicit and semi-implicit Euler methods for a
   blow-up problem in several space dimensions},
   journal={Computing},
   volume={76},
   date={2006},
   number={3-4},
   pages={325--352},
   issn={0010-485X},
   doi={10.1007/s00607-005-0136-0},
}
}

{\color{black}
\bib{JenseS2013}{article}{
   author={Jensen, Max},
   author={Smears, Iain},
   title={On the convergence of finite element methods for Hamilton-Jacobi-Bellman equations},
   journal={SIAM J. Numer. Anal.},
   volume={51},
   date={2013},
   number={1},
   pages={137--162},
   issn={0036-1429},
   doi={10.1137/110856198},
}
}

{\color{black}
\bib{KashiK2020a}{article}{
   author={Kashiwabara, Takahito},
   author={Kemmochi, Tomoya},
   title={Pointwise error estimates of linear finite element method for
   Neumann boundary value problems in a smooth domain},
   journal={Numer. Math.},
   volume={144},
   date={2020},
   number={3},
   pages={553--584},
   issn={0029-599X},
   doi={10.1007/s00211-019-01098-8},
}
}

{\color{black}
\bib{KashiK2020}{article}{
   author={Kashiwabara, Takahito},
   author={Kemmochi, Tomoya},
   title={Stability, analyticity, and maximal regularity for parabolic
   finite element problems on smooth domains},
   journal={Math. Comp.},
   volume={89},
   date={2020},
   number={324},
   pages={1647--1679},
   issn={0025-5718},
   doi={10.1090/mcom/3500},
}
}

{\color{black}
\bib{Knobl2006}{article}{
   author={Knobloch, Petr},
   title={Numerical solution of convection-diffusion equations using
   upwinding techniques satisfying the discrete maximum principle},
   conference={
      title={Proceedings of Czech-Japanese Seminar in Applied Mathematics
      2005},
   },
   book={
      series={COE Lect. Note},
      volume={3},
      publisher={Kyushu Univ. The 21 Century COE Program, Fukuoka},
   },
   date={2006},
   pages={69--76},
}
}

{\color{black}
\bib{LadySU1968}{book}{
    author={Lady{\v{z}}enskaja, O. A.},
    author={Solonnikov, V. A.},
    author={Ural{$'$}ceva, N. N.},
    title={Linear and quasilinear equations of parabolic type},
    series={Translated from the Russian by S. Smith. Translations of Mathematical Monographs, Vol. 23},
    publisher={American Mathematical Society, Providence, R.I.},
    date={1968},
    pages={xi+648},
}
}

{\color{black}
\bib{LeykeL2021}{article}{
   author={Leykekhman, Dmitriy},
   author={Li, Buyang},
   title={Weak discrete maximum principle of finite element methods in
   convex polyhedra},
   journal={Math. Comp.},
   volume={90},
   date={2021},
   number={327},
   pages={1--18},
   issn={0025-5718},
   doi={10.1090/mcom/3560},
}
}

{\color{black}
\bib{LiQXY2024}{article}{
   author={Li, Buyang},
   author={Qiu, Weifeng},
   author={Xie, Yupei},
   author={Yu, Wenshan},
   title={Weak discrete maximum principle of isoparametric finite element
   methods in curvilinear polyhedra},
   journal={Math. Comp.},
   volume={93},
   date={2024},
   number={345},
   pages={1--34},
   issn={0025-5718},
   doi={10.1090/mcom/3876},
}
}

{\color{black}
\bib{MizogS2023}{article}{
   author={Mizoguchi, Noriko},
   author={Souplet, Philippe},
   title={Singularity formation and regularization at multiple times in the
   viscous Hamilton-Jacobi equation},
   journal={Asymptot. Anal.},
   volume={133},
   date={2023},
   number={3},
   pages={291--353},
   issn={0921-7134},
   doi={10.3233/asy-221813},
}
}

{\color{black}
\bib{MizogS2023a}{article}{
   author={Mizoguchi, Noriko},
   author={Souplet, Philippe},
   title={Complete classification of gradient blow-up and recovery of
   boundary condition for the viscous Hamilton-Jacobi equation},
   journal={J. Geom. Anal.},
   volume={33},
   date={2023},
   number={2},
   pages={Paper No. 42, 147},
   issn={1050-6926},
   doi={10.1007/s12220-022-01002-1},
}
}

{\color{black}
\bib{Nakag1975}{article}{
    author={Nakagawa, Tomoyasu},
    title={Blowing up of a finite difference solution to {$u\sb{t}=u\sb{xx}+u^2$}},
    journal={Appl. Math. Optim.},
    volume={2},
    date={1975/76},
    number={4},
    pages={337--350},
    doi={10.1007/BF01448176},
    url={https://doi.org/10.1007/BF01448176},
    issn={0095-4616,1432-0606},
}
}

{\color{black}
\bib{Niren1959}{article}{
   author={Nirenberg, L.},
   title={On elliptic partial differential equations},
   journal={Ann. Scuola Norm. Sup. Pisa Cl. Sci. (3)},
   volume={13},
   date={1959},
   pages={115--162},
   issn={0391-173X},
}
}

{\color{black}
\bib{PorreS2017a}{article}{
    author={Porretta, Alessio},
    author={Souplet, Philippe},
    title={Analysis of the loss of boundary conditions for the diffusive {H}amilton-{J}acobi equation},
    journal={Ann. Inst. H. Poincar{\'e} Anal. Non Lin{\'e}aire},
    volume={34},
    date={2017},
    number={7},
    pages={1913--1923},
    doi={10.1016/j.anihpc.2017.02.001},
    url={https://doi.org/10.1016/j.anihpc.2017.02.001},
    issn={0294-1449},
}
}

{\color{black}
\bib{PorreS2017}{article}{
    author={Porretta, Alessio},
    author={Souplet, Philippe},
    title={The profile of boundary gradient blowup for the diffusive {H}amilton-{J}acobi equation},
    journal={Int. Math. Res. Not. IMRN},
    date={2017},
    number={17},
    pages={5260--5301},
    url={https://doi.org/10.1093/imrn/rnw154},
    issn={1073-7928},
}
}

{\color{black}
\bib{PorreS2020}{article}{
    author={Porretta, Alessio},
    author={Souplet, Philippe},
    title={Blow-up and regularization rates, loss and recovery of boundary conditions for the superquadratic viscous {H}amilton-{J}acobi equation},
    journal={J. Math. Pures Appl. (9)},
    volume={133},
    date={2020},
    pages={66--117},
    doi={10.1016/j.matpur.2019.02.014},
    url={https://doi.org/10.1016/j.matpur.2019.02.014},
    issn={0021-7824,1776-3371},
}
}

{\color{black}
\bib{RannaS1982}{article}{
   author={Rannacher, Rolf},
   author={Scott, Ridgway},
   title={Some optimal error estimates for piecewise linear finite element
   approximations},
   journal={Math. Comp.},
   volume={38},
   date={1982},
   number={158},
   pages={437--445},
   issn={0025-5718},
   doi={10.2307/2007280},
}
}

{\color{black}
\bib{Schat1980}{article}{
   author={Schatz, Alfred H.},
   title={A weak discrete maximum principle and stability of the finite
   element method in $L\sb{\infty }$\ on plane polygonal domains. I},
   journal={Math. Comp.},
   volume={34},
   date={1980},
   number={149},
   pages={77--91},
   issn={0025-5718},
   doi={10.2307/2006221},
}
}

{\color{black}
\bib{Soupl2002a}{article}{
    author={Souplet, Philippe},
    title={Gradient blow-up for multidimensional nonlinear parabolic equations with general boundary conditions},
    journal={Differential Integral Equations},
    volume={15},
    date={2002},
    number={2},
    pages={237--256},
    issn={0893-4983},
}
}

{\color{black}
\bib{SouplZ2006}{article}{
    author={Souplet, Philippe},
    author={Zhang, Qi S.},
    title={Global solutions of inhomogeneous {H}amilton-{J}acobi equations},
    journal={J. Anal. Math.},
    volume={99},
    date={2006},
    pages={355--396},
    doi={10.1007/BF02789452},
    issn={0021-7670},
}
}

{\color{black}
\bib{Tabat1977}{article}{
    author={Tabata, Masahisa},
    title={A finite element approximation corresponding to the upwind finite differencing},
    journal={Mem. Numer. Math.},
    date={1977},
    number={4},
    pages={47--63},
}
}

{\color{black}
\bib{Tabat1978}{article}{
   author={Tabata, Masahisa},
   title={Uniform convergence of the upwind finite element approximation for
   semilinear parabolic problems},
   journal={J. Math. Kyoto Univ.},
   volume={18},
   date={1978},
   number={2},
   pages={327--351},
   issn={0023-608X},
   doi={10.1215/kjm/1250522577},
}
}

{\color{black}
\bib{Thome2006}{book}{
   author={Thom\'ee, Vidar},
   title={Galerkin finite element methods for parabolic problems},
   series={Springer Series in Computational Mathematics},
   volume={25},
   edition={2},
   publisher={Springer-Verlag, Berlin},
   date={2006},
   pages={xii+370},
   isbn={978-3-540-33121-6},
   isbn={3-540-33121-2},
}
}

{\color{black}
\bib{ThomeW2008}{article}{
    author={Thom{\'e}e, Vidar},
    author={Wahlbin, Lars B.},
    title={On the existence of maximum principles in parabolic finite element equations},
    journal={Math. Comp.},
    volume={77},
    date={2008},
    number={261},
    pages={11--19},
    url={https://doi.org/10.1090/S0025-5718-07-02021-2},
    issn={0025-5718},
}
}

{\color{black}
\bib{Ushij2000}{article}{
   author={Ushijima, Takeo K.},
   title={On the approximation of blow-up time for solutions of nonlinear
   parabolic equations},
   journal={Publ. Res. Inst. Math. Sci.},
   volume={36},
   date={2000},
   number={5},
   pages={613--640},
   issn={0034-5318},
   doi={10.2977/prims/1195142812},
}
}

{\color{black}
\bib{MR1149111}{article}{
   author={van der Vorst, H. A.},
   title={Bi-CGSTAB: a fast and smoothly converging variant of Bi-CG for the
   solution of nonsymmetric linear systems},
   journal={SIAM J. Sci. Statist. Comput.},
   volume={13},
   date={1992},
   number={2},
   pages={631--644},
   issn={0196-5204},
   doi={10.1137/0913035},
}
}

{\color{black}
\bib{VandeHZG2010}{article}{
   author={VanderZee, Evan},
   author={Hirani, Anil N.},
   author={Zharnitsky, Vadim},
   author={Guoy, Damrong},
   title={A dihedral acute triangulation of the cube},
   journal={Comput. Geom.},
   volume={43},
   date={2010},
   number={5},
   pages={445--452},
   issn={0925-7721},
   doi={10.1016/j.comgeo.2009.09.001},
}
}

{\color{black}
\bib{ZhouS2017}{article}{
    author={Zhou, Guanyu},
    author={Saito, Norikazu},
    title={Finite volume methods for a {K}eller-{S}egel system: discrete energy, error estimates and numerical blow-up analysis},
    journal={Numer. Math.},
    volume={135},
    date={2017},
    number={1},
    pages={265--311},
    url={https://doi.org/10.1007/s00211-016-0793-2},
    issn={0029-599X},
}
}

\end{biblist}
\end{bibdiv}
\end{document}